\documentclass[10pt]{amsart}
\usepackage[margin=2.5cm]{geometry}
\usepackage{tikz-cd}
\usepackage{graphicx}					
\usepackage{amssymb}
\usepackage{amsmath}
\usepackage{subfig, graphicx}
\usepackage{mathrsfs}
\usepackage{amssymb, amsthm, indentfirst}
\usepackage{relsize}
\usepackage{hyperref}
\usepackage{stmaryrd}
\usepackage{lineno}

\DeclareMathOperator{\lcm}{lcm}

\numberwithin{equation}{section}
\usepackage{color}
\theoremstyle{plain}
\newtheorem{thm}{Theorem}[section]

\newtheorem{conj}[thm]{Conjecture}
\newtheorem{lemma}[thm]{Lemma}
\newtheorem{prop}[thm]{Proposition}

\theoremstyle{definition}

\newtheorem{rmk}[thm]{Remark}
\def\Gal{\operatorname{Gal}}

\newtheorem*{hypothesis*}{Hypothesis}

\author{SHIH-YU CHEN}
\address{Institute of Mathematics~\\Academia Sinica~\\ 6F, Astronomy-Mathematics Building, No.\,1, Sec.\,4, Roosevelt Road, Taipei 10617, Taiwan, ROC}
\email{sychen0626@gate.sinica.edu.tw}

\def\SL{{\rm{SL}}}
\def\GL{{\rm{GL}}}
\def\GSp{{\rm GSp}}

\def\Sp{{\rm Sp}}

\def\o{\frak{o}}

\def\A{{\mathbb A}}
\def\C{{\mathbb C}}
\def\E{{\mathbb E}}
\def\F{{\mathbb F}}
\def\K{{\mathbb K}}
\def\L{{\mathbb L}}
\def\R{{\mathbb R}}
\def\Q{{\mathbb Q}}
\def\Z{{\mathbb Z}}

\def\<{\langle}
\def\>{\rangle}
\def\G{\mathbf{G}}

\def\bp{\begin{pmatrix}}
\def\ep{\end{pmatrix}}

\def\<{\langle}
\def\>{\rangle}

\def\GL{\operatorname{GL}}

\def\GSp{\operatorname{GSp}}

\def\SL{\operatorname{SL}}

\def\Sp{\operatorname{Sp}}

\def\KK{\mathbb{K}}

\def\1{\mathbf{1}}

\def\itPi{\mathit{\Pi}}

\def\itSigma{\mathit{\Sigma}}
\def\itPhi{\Phi}

\title{Algebraicity of critical values of triple product $L$-functions in the balanced case}
\begin{document}

\begin{abstract}
The algebraicity of critical values of triple product $L$-functions in the balanced case was proved by Garrett and Harris, under the assumption that the critical points are on the right and away from center of the critical strip. The missing right-half critical points correspond to certain holomorphic Eisenstein series outside the range of absolute convergence. The remaining difficulties are construction of these holomorphic Eisenstein series and verification of the non-vanishing of the corresponding non-archimedean local zeta integrals.
In this paper, we address these problems and complement the result of Garrett and Harris to all critical points.
As a consequence, we obtain new cases of Deligne's conjecture for symmetric cube $L$-functions of Hilbert modular forms.
\end{abstract}

\maketitle

\section{Introduction}

The study of special values of automorphic $L$-functions is an important subject in number theory.
In \cite{Garrett1987}, Garrett discovered an integral representation of the triple product $L$-function attached to a triplet $(f_1,f_2,f_3)$ of normalized elliptic newforms. 
Subsequent works of Piatetski-Shapiro--Rallis \cite{PSR1987}, Ikeda \cite{Ikeda1989}, \cite{Ikeda1992}, \cite{Ikeda1998}, \cite{Ikeda1999}, and Ramakrishnan \cite{Rama2000} studied the analytic properties of the triple product $L$-functions both globally and locally.
More precisely, for $i=1,2,3$, let $f_i$ be a normalized elliptic newform of weight $\kappa_i$, level $N_i$, and nebentypus $\omega_i$. 
We have the triple product $L$-function associated to the triplet $(f_1,f_2,f_3)$ defined by an Euler product
\[
L(s,f_1\times f_2 \times f_3) = \prod_p L_p(s,f_1\times f_2 \times f_2).
\]
The Euler product convergent absolutely for ${\rm Re}(s)$ sufficiently large, admits meromorphic continuation to $s \in \C$, and satisfies a functional equation relating $L(s,f_1\times f_2\times f_3)$ and $L(\kappa_1+\kappa_2+\kappa_3-2-s,f_1^\vee\times f_2^\vee \times f_3^\vee)$.
Here $f_i^\vee$ is the normalized elliptic newform dual to $f_i$ under the Hecke action. 
When $\kappa_1=\kappa_2=\kappa_3=\kappa$ and $N_1=N_2=N_3=1$, Garrett proved the algebraicity of the rightmost critical value $L(2\kappa-2,f_1\times f_2 \times f_3)$ in terms of the product of Petersson norms of $f_1$, $f_2$, and $f_3$.
In general, we consider the algebraicity of critical values when the triplet of weights $(\kappa_1,\kappa_2,\kappa_3)$ is in the balanced range, that is, 
\[
\kappa_1+\kappa_2+\kappa_3 >2\max\{\kappa_1,\kappa_2,\kappa_3\}.
\]
In the balanced case, the conjectural Deligne's period attached to the triple product $L$-function was computed by Blasius \cite{Blasius1987} and is equal to the product of Petersson norms of $f_1$, $f_2$, and $f_3$ times some elementary factors. Under the same assumptions as in \cite{Garrett1987}, Orloff \cite{Orloff1987} and Satoh \cite{Satoh1987} proved the algebraicity of all critical values. 
The triple product $L$-function has a central critical point $m_0 = \tfrac{\kappa_1+\kappa_2+\kappa_3}{2}-1$ when $\kappa_1+\kappa_2+\kappa_3$ is even.
In this case, Harris and Kudla \cite{HK1991} proved the algebraicity of the central critical value under the assumption that $\omega_1\omega_2\omega_3$ is trivial.
In \cite{GH1993}, Garrett and Harris proved the algebraicity in general case over totally real number fields for right-half critical points away from center of the critical strip.
To be precise, they proved the algebraicity of $L(m,f_1\times f_2 \times f_3)$ when $m \geq m_0$ and $m \notin \{m_0,m_0+\tfrac{1}{2},m_0+1\}$. B\"ocherer and Schulze-Pillot \cite{BSP1996} proved the algebraicity of all critical values under the assumptions that $\lcm(N_1,N_2,N_3)>1$ is square-free and $\omega_1=\omega_2=\omega_3=1$.
In their construction of $p$-adic triple product $L$-function, Hsieh and Yamana \cite{HY2019} obtain the algebraicity of all critical values as a byproduct under the assumption that $f_1$, $f_2$, and $f_3$ are simultaneously ordinary at an odd prime $p$. 
For numerical computation of critical values of triple product $L$-functions, there are results of Mizumoto \cite{Mizumoto2000} and Ibukiyama--Katsurada--Poor--Yuen \cite{IKPY2014}. For explicit central value formula for the triple product $L$-function, we refer to the results of Gross and Kudla \cite{GK1992}, \cite[\S\,5]{BSP1996}, and the author and Cheng \cite[\S\,6]{CC2019}.

The purpose of this paper is to complement the theorem of Garrett and Harris in \cite{GH1993}. As pointed out in \cite[Remarks 3.4.8 and 4.6.3]{GH1993}, the remaining difficulties for the right-half critical pointes $m \in \{m_0,m_0+\tfrac{1}{2},m_0+1\}$ are construction of holomorphic Siegel Eisenstein series outside the range of absolute convergence and verification of the non-vanishing of the corresponding non-archimedean local zeta integral.
These issues are addressed in the present paper. 
We are inspired by the construction of Eisenstein series due to Shimura \cite{Shimura1997}.
More precisely, we choose an auxiliary good place at which the local section defining the Eisenstein series is supported in the big cell.
For the local computation of zeta integral, we follow the method employed in \cite{HY2019}. It turns out that when $m \neq m_0+\tfrac{1}{2}$, we can find input data so that the resulting local zeta integral is non-vanishing at $s=m$. Whereas when $m=m_0+\tfrac{1}{2}$ and either $(\omega_1\omega_2\omega_3)^2 \neq 1$ or $f_1,f_2,f_3$ are not simultaneously CM by an imaginary quadratic field, we need the Sato--Tate conjecture to guarantee the existence of good place such that some local zeta integral is non-vanishing at $s=m$.
For $m = m_0 +\tfrac{1}{2}$ and $f_1,f_2,f_3$ are simultaneously CM by an imaginary quadratic field, we can reduce to the previous case $(\omega_1\omega_2\omega_3)^2 \neq 1$ by showing that the algebraicity holds for $L(m,f_1\times f_2 \times f_3)$ if and only if it holds for $L(m,f_1\times f_2 \times f_3 \otimes\chi)$ for any even Dirichlet character $\chi$. The algebraicity for the left-half critical points then follows from the functional equation of triple product $L$-function.
In this paper, we also consider the twisted cases as in \cite{GH1993}. We consider the triple product $L$-function associated to a Hilbert cusp newform over a totally real \'etale cubic algebra $\E/\F$ for some totally real number field $\F$. For the non-split case $\E \neq \F\times\F\times\F$, the functional equation for the triple product $L$-function follows from \cite{PSR1987} and our previous result \cite{Chen2021b} on gamma factors. 

The main motivation behind this paper is application to Deligne's conjecture for symmetric power $L$-functions of Hilbert modular forms (cf.\,\cite[\S\,3]{RS2007c}).
As proved in \cite[\S\,6]{GH1993}, Deligne's conjecture for symmetric cube $L$-functions holds for critical points away from center of the critical strip. As a consequence, we obtain new cases for symmetric cube $L$-functions. Moreover, combined with the results of Grobner--Raghuram \cite{GR2014} and Harder--Raghuram \cite{HR2020}, we also obtain the algebraicity of the (possible) central value of the symmetric cube $L$-functions.
In \cite{Morimoto2021}, the result of Garrett and Harris was used by Morimoto as one of the key ingredients to prove the algebraicity of critical values of symmetric fourth $L$-functions of Hilbert modular forms. 
Because we can now prove the complete result concerning \cite{GH1993}, the condition on the weight of Hilbert modular forms can be further relaxed. 

\subsection{Main results}
\subsubsection{Triple product $L$-functions}
Let $\F$ be a totally real number field with $d=[\F:\Q]$ and $\E$ a totally real \'etale cubic algebraic over $\F$.
Denote by $S_\infty$ the set of archimedean places of $\F$.
Let $\itPi = \bigotimes_v \itPi_v$ be a cohomological irreducible cuspidal automorphic representation of $({\rm R}_{\E/\F}\GL_2)(\A_\F)=\GL_2(\A_\E)$ with central character $\omega_\itPi$.
Put $\omega = \omega_\itPi\vert_{\A_\F^\times}$.
Let $f_{\itPi}$ and $f_{\itPi^\vee}$ be the normalized newforms of $\itPi$ and $\itPi^\vee$.  
The Petersson norm of $f_{\itPi}$ is defined by
\begin{align}\label{E:Petersson}
\Vert f_{\itPi} \Vert = \int_{\A_\E^\times\GL_2(\E)\backslash \GL_2(\A_\E)} f_{\itPi}(g)f_{\itPi^\vee}(g\cdot{\rm diag}(-1,1)_\infty)\,dg^{\rm Tam}.
\end{align}
Here $dg^{\rm Tam}$ is the Tamagawa measure on $\A_\E^\times \backslash \GL_2(\A_\E)$.
Let 
\[
L(s,\itPi, {\rm As}) = \prod_{v}L(s,\itPi,{\rm As})
\]
be the triple product $L$-function of $\itPi$ defined by the Asai cube representation (see \S\,\ref{S:triple L-function} for the precise definition)
\[
{\rm As} : {}^L({\rm R}_{\E/\F}\GL_2) \longrightarrow \GL(\C^2\otimes\C^2\otimes\C^2).
\]
The Euler product is convergent absolutely for ${\rm Re}(s)$ sufficiently large and admits meromorphic continuation to $s \in \C$.
We denote by $L^{(\infty)}(s,\itPi,{\rm As})$ the $L$-function obtained by excluding the archimedean local factors. A critical point for $L(s,\itPi,{\rm As})$ is a half-integer $m+\tfrac{1}{2} \in \Z+\tfrac{1}{2}$ which is not a pole of the archimedean local factors $L(s,\itPi_v,{\rm As})$ and $L(1-s,\itPi_v^\vee,{\rm As})$ for all $v \in S_\infty$. 
For $v \in S_\infty$, we have
\[
\itPi_v = \boxtimes_{i=1}^3 D_{\kappa_{i,v}} \otimes |\mbox{ }|^{{\sf w}_{i,v}/2}
\]
as representations of $({\rm R}_{\E_v/\F_v}\GL_2)(\F_v) = \GL_2(\R)^3$ for some $\kappa_{i,v} \in \Z_{\geq 2}$ and ${\sf w}_{i,v} \in \Z$ such that $\kappa_{i,v} \equiv {\sf w}_{i,v} \,({\rm mod}\,2)$. Here $D_\kappa$ is the discrete series representation of $\GL_2(\R)$ with weight $\kappa \in \Z_{\geq 2}$.
We assume the weights satisfy the balanced condition:
\begin{align}\label{E:balanced}
\kappa_{1,v}+\kappa_{2,v}+\kappa_{3,v}>2\max\{\kappa_{1,v},\kappa_{2,v},\kappa_{3,v}\}
\end{align}
for all $v \in S_\infty$.
In this case, $m+\tfrac{1}{2}$ is critical if and only if
\[
- \left( \tfrac{\kappa_{1,v}+\kappa_{2,v}+\kappa_{3,v}-2\max\{\kappa_{1,v},\kappa_{2,v},\kappa_{3,v}\}}{2} \right) + 1 \leq m+\tfrac{{\sf w}}{2} \leq \tfrac{\kappa_{1,v}+\kappa_{2,v}+\kappa_{3,v}-2\max\{\kappa_{1,v},\kappa_{2,v},\kappa_{3,v}\}}{2} -1
\]
for all $v \in S_\infty$. Here ${\sf w} \in \Z$ is the integer such that $|\omega| = |\mbox{ }|_{\A_\F}^{\sf w}$.
It is clear that a necessary and sufficient condition for the existence of critical point is that
\begin{align}\label{E:balanced 2}
\kappa_{1,v}+\kappa_{2,v}+\kappa_{3,v}\geq 2\max\{\kappa_{1,v},\kappa_{2,v},\kappa_{3,v}\}+2
\end{align}
for all $v \in S_\infty$.
Note that the triple product $L$-function has a central critical point $s = \tfrac{1-{\sf w}}{2}$ when ${\sf w}$ is even.
Also one can deduce from the result \cite[Theorem 2.7]{Ikeda1992} of Ikeda that the triple product $L$-function is holomorphic at critical points.
For $\sigma \in {\rm Aut}(\C)$, let ${}^\sigma\!\itPi$ be the irreducible admissible representation of $({\rm R}_{\E/\F}\GL_2)(\A_\F)$ defined by
\[
{}^\sigma\!\itPi = {}^\sigma\!\itPi_{\infty} \otimes {}^\sigma\!\itPi_f,
\]
where ${}^\sigma\!\itPi_f$ is the $\sigma$-conjugate of $\itPi_f = \bigotimes_{v\nmid \infty} \itPi_v$ and ${}^\sigma\!\itPi_{\infty}$ is the representation of $({\rm R}_{\E/\F}\GL_2)(\F_\infty)$ so that its $v$-component is equal to $\itPi_{\sigma^{-1} \circ v}$ for $v \in S_\infty$.
Since $\itPi$ is cohomological, it is known that ${}^\sigma\!\itPi$ is cuspidal automorphic.
In \cite{GH1993}, Garrett and Harris proved the following result on the algebraicity of critical values of the triple product $L$-function.

\begin{thm}[Garrett--Harris]\label{T:GH}
Assume the balanced condition (\ref{E:balanced 2}) is satisfied.
Let $m+\tfrac{1}{2} \in \Z+\tfrac{1}{2}$ be a critical point for $L(s,\itPi,{\rm As})$ such that $m+\tfrac{{\sf w}}{2} \geq \tfrac{3}{2}$.
Then we have
\begin{align*}
&\sigma \left( \frac{L^{(\infty)}(m+\tfrac{1}{2},\itPi,{\rm As})}{(|D_\E|/|D_\F|)^{1/2}\cdot(2\pi\sqrt{-1})^{4dm+\sum_{v \in S_\infty}(\kappa_{1,v}+\kappa_{2,v}+\kappa_{3,v}+2{\sf w}+2)}\cdot(\sqrt{-1})^{d{\sf w}}\cdot G(\omega)^2\cdot \Vert f_{\itPi} \Vert}\right)\\
&=\frac{L^{(\infty)}(m+\tfrac{1}{2},{}^\sigma\!\itPi,{\rm As})}{(|D_\E|/|D_\F|)^{1/2}\cdot(2\pi\sqrt{-1})^{4dm+\sum_{v \in S_\infty}(\kappa_{1,v}+\kappa_{2,v}+\kappa_{3,v}+2{\sf w}+2)}\cdot(\sqrt{-1})^{d{\sf w}}\cdot G({}^\sigma\!\omega)^2\cdot \Vert f_{{}^\sigma\!\itPi} \Vert}
\end{align*}
for all $\sigma \in {\rm Aut}(\C)$. Here $|D_\F|$ and $|D_\E|$ are the absolute discriminant of $\E/\Q$ and $\F/\Q$, respectively, and $G(\omega)$ is the Gauss sum of $\omega = \omega_\itPi \vert_{\A_\F^\times}$.
\end{thm}

Following is the main result of this paper.
We extend the theorem of Garrett and Harris to the remaining critical $L$-values.

\begin{thm}\label{T:main 3}
Theorem \ref{T:GH} holds for all critical points.
\end{thm}

\subsubsection{Symmetric cube $L$-functions}\label{SS:symmetric cube}

Let $\itPi$ be a cohomological irreducible cuspidal automorphic representation of $\GL_2(\A_\F)$ with central character $\omega_\itPi$. Let ${\sf w}$ be the integer such that $|\omega_\itPi| = |\mbox{ }|_{\A_\F}^{\sf w}$. For $v \in S_\infty$, we have
\[
\itPi_v = D_{\kappa_{v}} \otimes |\mbox{ }|^{{\sf w}/2}
\]
for some $\kappa_v \geq 2$ with $\kappa_v \equiv {\sf w}\,({\rm mod}\,2)$.
Let $f_\itPi$ be the normalized newform of $\itPi$ and $\Vert f_\itPi \Vert$ its Petersson norm defined as in (\ref{E:Petersson}).
For a finite order Hecke character $\chi$ of $\A_\F^\times$, let 
\[
L(s,\itPi,{\rm Sym}^3 \otimes \chi)
\]
be the twisted symmetric cube $L$-function of $\itPi$ by $\chi$. We denote by $L^{(\infty)}(s,\itPi,{\rm Sym}^3 \otimes \chi)$ the $L$-function obtained by excluding the archimedean $L$-factors. 
The set of critical points for $L(s,\itPi,{\rm Sym}^3 \otimes \chi)$ is given by
\begin{align*}
\left\{ m+\tfrac{1}{2} \in \Z+\tfrac{1}{2}\,\left\vert\,\tfrac{-\min_{v \in S_\infty}\{\kappa_v\}-3{\sf w}}{2}+1 \leq m \leq \tfrac{\min_{v \in S_\infty}\{\kappa_v\}-3{\sf w}}{2}-1\right\}\right..
\end{align*}

For $\sigma \in {\rm Aut}(\C)$ and $\underline{\varepsilon} \in \{\pm1\}^{S_\infty}$, let $p({}^\sigma\!\itPi,\underline{\varepsilon}) \in \C^\times$ be the non-zero complex number defined in \cite[Theorem 4.3]{Shimura1978} (cf.\,Remark \ref{R:Shimura}).
We have the following conjecture proposed by Deligne \cite[\S\,7]{Deligne1979} on the algebraicity of the critical values of $L(s,\itPi,{\rm Sym}^3 \otimes \chi)$ in terms of $p(\itPi,\underline{\varepsilon})$ and $\Vert f_\itPi \Vert$.

\begin{conj}[Deligne]\label{C:DC}
Let $m+\tfrac{1}{2}$ be a critical point for $L(s,\itPi,{\rm Sym}^3 \otimes \chi)$. Then we have
\begin{align*}
&\sigma\left(\frac{L^{(\infty)}(m+\tfrac{1}{2},\itPi,{\rm Sym}^3 \otimes \chi)}{(2\pi\sqrt{-1})^{2dm+\sum_{v\in S_\infty}(2\kappa_v+3{\sf w})}\cdot(\sqrt{-1})^{d{\sf w}}\cdot G(\chi\omega_\itPi)^2\cdot p(\itPi,(-1)^{m+1}{\rm sgn}(\chi))^2\cdot\Vert f_\itPi \Vert}\right) \\
&= \frac{L^{(\infty)}(m+\tfrac{1}{2},{}^\sigma\!\itPi,{\rm Sym}^3 \otimes {}^\sigma\!\chi)}{(2\pi\sqrt{-1})^{2dm+\sum_{v\in S_\infty}(2\kappa_v+3{\sf w})}\cdot(\sqrt{-1})^{d{\sf w}}\cdot G({}^\sigma\!\chi{}^\sigma\!\omega_\itPi)^2\cdot p({}^\sigma\!\itPi,(-1)^{m+1}{\rm sgn}(\chi))^2\cdot\Vert f_{{}^\sigma\!\itPi} \Vert}
\end{align*}
for all $\sigma \in {\rm Aut}(\C)$. Here ${\rm sgn}(\chi) \in \{\pm1\}^{S_\infty}$ is the signature of $\chi$. 
\end{conj}

\begin{rmk}
The conjecture is true when $\itPi$ is of CM-type (cf.\,\cite[\S\,4]{RS2007c}).
\end{rmk}

As a consequence of Theorem \ref{T:GH}, Garrett and Harris proved the conjecture when the critical points are on the right and away from center of the critical strip. More precisely, we have the following result \cite[Theorem 6.2]{GH1993}.
\begin{thm}[Garrett--Harris]\label{T:GH DC}
If $\kappa_v \geq 5$ for all $v \in S_\infty$, then Conjecture \ref{C:DC} holds for critical points $m+\tfrac{1}{2}$ such that $m+\tfrac{3{\sf w}}{2}\geq \tfrac{3}{2}$.
\end{thm}

Following is the second main result of this paper.
Based on Theorem \ref{T:main 3}, we prove the conjecture when $\min_{v \in S_\infty}\{\kappa_v\}\geq 3$. For the (possible) central critical point $s=\tfrac{1-3{\sf w}}{2}$, we also use the results of Grobner--Raghuram \cite{GR2014} and Harder--Raghuram \cite{HR2020} to deduce the algebraicity.

\begin{thm}\label{T:DC}
If $\kappa_v \geq 3$ for all $v \in S_\infty$, then Conjecture \ref{C:DC} holds for all critical points.
\end{thm}

\begin{rmk}
When $\min_{v \in S_\infty}\{\kappa_v\}=2$, there is only one critical point $s=\tfrac{1-3{\sf w}}{2}$ and Conjecture \ref{C:DC} holds for all $\chi$ such that the central critical value $L(\tfrac{1-{\sf w}}{2},\itPi\otimes\chi)$ is non-zero. 
\end{rmk}

This paper is organized as follows. In \S\,\ref{S:Eisenstein series}, we construct holomorphic Siegel Eisenstein series on $\GSp_{2n}(\A_\F)$ outside the region of absolute convergence. With certain local assumption, we extend the result of Harris on the algebraicity of Siegel Eisenstein series. In \S\,\ref{S:triple}, we recall the definition of triple product $L$-function and Garrett's integral representation from adelic point of view.
We prove our main results in \S\,\ref{S:proofs}. In \S\,\ref{SS:proof 1}, we firstly explain how to reduce the problem to the construction in \S\,\ref{S:Eisenstein series} for $\GSp_6$ and computation of non-archimedean local zeta integrals. Subject to Proposition \ref{P:local zeta} on certain explicit local zeta integrals, we then deduce Theorem \ref{T:main 3}. 
The case when $m+\tfrac{\sf w}{2} = \tfrac{1}{2}$ is particularly involved as the local zeta integral in Proposition \ref{P:local zeta} might be zero in this case. We treat this case in \S\,\ref{SS:Case 1}-\ref{SS:Case 3}.
In \S\,\ref{S:local zeta integral}, we compute local zeta integrals and prove Proposition \ref{P:local zeta}.

\subsection{Notation}\label{SS:notation}
Fix a totally real number field $\F$ with $[\F:\Q]=d$.
Let $D_\F$ be the discriminant of $\F$.
Let $\A_\F$ be the ring of adeles of $\F$ and $\A_{\F,f}$ be its finite part. 
Let $\hat{\o}_\F$ be the maximal compact subring of $\A_{\F,f}$.
Let $\psi_\Q=\bigotimes_v\psi_{v}$ be the standard additive character of $\Q\backslash \A_\Q$ defined so that
\begin{align*}
\psi_{p}(x) & = e^{-2\pi \sqrt{-1}\,x} \mbox{ for }x \in \Z[p^{-1}],\\
\psi_{\infty}(x) & = e^{2\pi \sqrt{-1}\,x} \mbox{ for }x \in \R.
\end{align*}
Let $\psi_\F = \psi_\Q\circ{\rm tr}_{\F/\Q}$ and call it the standard additive character of $\F \backslash \A_\F$.
Let $S_\infty$ be the set of archimedean places of $\F$. 
For $v \in S_\infty$, let $\iota_v$ be the real embedding of $\F$ associated to $v$ and identify $\F_v$ with $\R$ via $\iota_v$.
Let $v$ be a place of $\F$. If $v$ is a finite place, let $\frak{o}_{\F_v}$, $\varpi_v$, and $q_v$ be the maximal compact subring of $\F_v$, a generator of the maximal ideal of $\frak{o}_{\F_v}$, and the cardinality of $\frak{o}_{\F_v} / \varpi_v\frak{o}_{\F_v}$. Let $|\mbox{ }|_v$ be the absolute value on $\F_v$ normalized so that $|\varpi_v|_v = q_v^{-1}$. If $v \in S_\infty$, let $|\mbox{ }|_v$ be the be the ordinary absolute value on $\R$.
Let $|\mbox{ }|_{\A_\F} = \prod_v|\mbox{ }|_v$ be the adelic norm on $\A_\F$.

Let $\chi$ be an algebraic Hecke character of $\A_\F^\times$.
The signature of $\chi$ at $v \in S_\infty$ is the value $\chi_v(-1) \in \{\pm1\}$. The signature ${\rm sgn}(\chi)$ of $\chi$ is the sequence of signs $(\chi_v(-1))_{v \in S_\infty}$.
We say $\chi$ has parallel signature if it has the same signature at all real places.
The Gauss sum $G(\chi)$ of $\chi$ is defined by
\[
G(\chi) = |D_\F|^{-1/2}\prod_{v \nmid \infty}\varepsilon(0,\chi_v,\psi_v),
\]
where $\psi_\F = \bigotimes_v\psi_v$ and $\varepsilon(s,\chi_v,\psi_v)$ is the $\varepsilon$-factor of $\chi_v$ with respect to $\psi_v$ defined in \cite{Tate1979}.
For $\sigma \in {\rm Aut}(\C)$, let ${}^\sigma\!\chi$ be the unique algebraic Hecke character  of $\A_\F^\times$ such that ${}^\sigma\!\chi(x) = \sigma(\chi(x))$ for $x \in \A_{\F,f}^\times$. Note that ${\rm sgn}(\chi) = {\rm sgn}({}^\sigma\!\chi)$.
It is easy to verify that 
\begin{align}\label{E:Galois Gauss sum}
\begin{split}
\sigma(G(\chi)) = {}^\sigma\!\chi(u_\sigma)G({}^\sigma\!\chi),\quad
\sigma\left(\frac{G(\chi\chi')}{G(\chi)G(\chi')}\right) = \frac{G({}^\sigma\!\chi{}^\sigma\!\chi')}{G({}^\sigma\!\chi)G({}^\sigma\!\chi')}
\end{split}
\end{align}
for algebraic Hecke characters $\chi,\chi'$ of $\A_\F^\times$,
where $u_\sigma \in \prod_p \Z_p^\times$ is the unique element such that $\sigma(\psi_{\Q}(x)) = \psi_{\Q}(u_\sigma x)$ for $x \in \A_{\Q,f}$.

Let $\sigma \in {\rm Aut}(\C)$. Define the $
\sigma$-linear action on $\C(X)$, which is the field of formal Laurent series in variable $X$ over $\C$, as follows:
\[
{}^\sigma\!P(X) = \sum_{n \gg -\infty}^\infty\sigma(a_n) X^n
\]
for $P(X) = \sum_{n \gg -\infty}^\infty a_n X^n \in \C(X)$.

\section{Holomorphic Siegel Eisenstein series}\label{S:Eisenstein series} 

In this section, we construct holomorphic Siegel Eisenstein series on $\GSp_{2n}(\A_\F)$ outside the range of absolute convergence. The main result of this section is Proposition \ref{P:holomorphy}, where we extend the result of Harris on algebraicity of Siegel Eisenstein series. We will apply the construction to $\GSp_6$ in \S\,\ref{S:triple} and \S\,\ref{S:proofs}.

\subsection{Hilbert--Siegel modular forms}\label{S:H-S modular form}

Let $\GSp_{2n}$ be the symplectic similitude group of degree $n$ defined by
\[
\GSp_{2n} = \left\{g \in \GL_{2n}\,\left\vert\,gJ_n{}^t\!g = \nu(g)J_n\right\}\right.,\quad J_n = \bp 0 & {\bf 1}_n \\ -{\bf 1}_n & 0\ep.
\]
For $\nu\in\GL_1$, $a \in \GL_n$, and $x,y \in {\rm Sym}_n$, define ${\bf m}(a,\nu),{\bf n}(x),{\bf n}^-(y) \in \GSp_{2n}$ by
\begin{align}\label{E:group elements}
{\bf m}(a,\nu) = \bp a & 0 \\ 0 & \nu {}^t\!a^{-1}\ep ,\quad {\bf n}(x) = \bp {\bf 1}_n& x \\ 0 & {\bf 1}_n \ep,\quad {\bf n}^-(y) = \bp {\bf 1}_n& 0 \\ y & {\bf 1}_n \ep.
\end{align}
We write ${\bf m}(a) = {\bf m}(a,1)$. When $n=1$, we also write ${\bf t}(a) = {\bf m}(a,a)$ for $a \in \GL_1$.
Let $P_n$ be the standard Siegel parabolic subgroup of $\GSp_{2n}$ defined by
\[
P_n =\left. \left\{ {\bf m}(a,\nu){\bf n}(x)\,\right\vert\, \nu \in\GL_1,\,a\in\GL_n,\,x\in{\rm Sym}_n\right\}.
\]
Let $\GSp_{2n}^+(\R)$ be the subgroup of $\GSp_{2n}(\R)$ consisting elements with positive similitude. 
We identify the compact unitary group ${\rm U}(n)$ with a maximal compact subgroup of $\Sp_{2n}(\R)$ via the map
\[
{\rm U}(n)\ni a+b{\sqrt{-1}} \longmapsto \bp a & b \\ -b & a \ep \in \Sp_{2n}(\R).
\]
For $v \in S_\infty$, we identify $\GSp_{2n}(\F_v)$ with $\GSp_{2n}(\R)$ through the real embedding $\iota_v$ of $\F$ associated to $v$.

Let $\ell \geq 1$ and $\chi$ be a Hecke character of $\A_\F^\times$.
An automorphic form $\varphi$ on $\GSp_{2n}(\A_\F)$ is holomorphic of weight $\ell$ and has central character $\chi$ if the following conditions are satisfied:
\begin{itemize}
\item $\varphi(gak_\infty) = \chi(a)(\det k_\infty)^\ell\varphi(g)$ for all $a \in \A_\F^\times$, $k_\infty \in {\rm U}(n)^{S_\infty}$, and $g \in \GSp_{2n}(\A_\F)$.
\item $(\frak{p}^-)^{S_\infty}\cdot\varphi = 0$, where
\[
\frak{p}^- = \left.\left\{ \bp A & -\sqrt{-1}\,A \\ -\sqrt{-1}\,A & -A\ep\,\right\vert\, A \in {\rm Sym}_n(\C) \right\} \subset {\rm Lie}(\GSp_{2n}(\R))_\C.
\]
\end{itemize}
Note that a necessary condition for $\varphi$ to be non-zero is that $\chi_v(-1) = (-1)^{n\ell}$ for all $v \in S_\infty$.
For an automorphic form $\varphi$ on $\GSp_{2n}(\A_\F)$, we have the Fourier expansion:
\[
\varphi = \sum_{B \in {\rm Sym}_n(\F)}W_B(\varphi),
\]
where $W_B(\varphi)$ is the $B$-th Fourier coefficient defined by
\begin{align*}
W_B(g,\varphi) = \int_{{\rm Sym}_n(\F)\backslash{\rm Sym}_n(\A_\F)}\varphi\left( {\bf n}(x) g\right)\overline{\psi_\F({\rm tr}(Bx))}\,dx^{\rm Tam}
\end{align*}
for $g \in \GSp_{2n}(\A_\F)$. 
Here $dx^{\rm Tam}$ is the Tamagawa measure on ${\rm Sym}_n(\A_\F)$.
Assume $\varphi$ is holomorphic of weight $\ell$ and has central character $\chi$, then we have
\begin{align}\label{E:q-expansion}
e^{2\pi\sum_{v \in S_\infty}{\rm tr}({}^t\!a_v\iota_v(B)a_v)}(\det a_\infty)^{-\ell}\cdot W_B\left({\bf m}(a_\infty)\cdot g_f,\varphi\right)
 = e^{2\pi {\rm tr}_{\F/\Q}({\rm tr}(B))}\cdot W_B(g_f,\varphi)
\end{align}
for all $g_f \in \GSp_{2n}(\A_{\F,f})$ and $a_\infty = (a_v)_{v \in S_\infty} \in \GL_n(\R)^{S_\infty}$.
Also $W_B(\varphi) \neq 0$ only when $B$ is totally positive semi-definite, that is, $\iota_v(B)$ is positive semi-definite for all $v \in S_\infty$.
We have the following result of Harris \cite[\S\,6]{Harris1986} (see also \cite[Appendix A.4]{GH1993} and \cite[(1.1.13)]{BHR1994}) on the algebraicity of holomorphic automorphic forms.

\begin{thm}[$q$-expansion principle]\label{T:q-expansion}
Let $\ell \geq 1$ and $\chi$ be an algebraic Hecke character of $\A_\F^\times$ with parallel signature $(-1)^{n\ell}$ and $|\chi| = |\mbox{ }|_{\A_\F}^{\sf u}$ for some ${\sf u}\equiv n\ell\,({\rm mod}\,2)$. 
Let $\varphi$ be a holomorphic automorphic form of weight $\ell$ and has central character $\chi$.
For $\sigma \in {\rm Aut}(\C)$, there is a holomorphic automorphic form ${}^\sigma\!\varphi$ of weight $\ell$ and has central character ${}^\sigma\!\chi$ uniquely determined by
\begin{align*}
\sigma\left( e^{2\pi {\rm tr}_{\F/\Q} ({\rm tr}(B))}\cdot W_{B} \left(\bp u_\sigma^{-1}{\bf 1}_n & 0 \\ 0 & {\bf 1}_n\ep g_f,\varphi \right)\right)= e^{2\pi {\rm tr}_{\F/\Q} ({\rm tr}(B))} \cdot W_{B} (g_f,{}^\sigma\!\varphi)
\end{align*}
for $B \in {\rm Sym}_n(\F)$ and $g_f \in \GSp_{2n}(\A_{\F,f})$.
Here $u_\sigma \in \widehat{\Z}^\times$ is the unique element such that $\sigma(\psi_{\Q}(x)) = \psi_{\Q}(u_\sigma x)$ for $x \in \A_{\Q,f}$.
\end{thm}

We recall the classical reformulation in terms of Hilbert--Siegel modular forms.
Let $\frak{H}_n$ be the Siegel upper-half space of degree $n$ defined by
\[
\frak{H}_n = \left\{ x+y\sqrt{-1} \in {\rm Sym}_n(\C) \,\vert\, y\mbox{ positive definite}\right\}.
\]
Let $\varphi$ be a holomorphic automorphic form on $\GSp_{2n}(\A_\F)$ of weight $\ell$ and has central character $\chi$. Let ${\bf f}_\varphi : \frak{H}_n^{S_\infty} \times \GSp_{2n}(\A_{\F,f}) \rightarrow \C$ be the Hilbert--Siegel modular form associated to $\varphi$ defined by
\begin{align}\label{E:H-S modular form}
{\bf f}_\varphi(x+y\sqrt{-1} ,g_f) = (\det y)^{-\ell/2}\cdot \varphi({\bf n}(x){\bf m}(y^{1/2})\cdot g_f).
\end{align}
Then the condition $(\frak{p}^-)^{S_\infty}\cdot\varphi=0$ is equivalent to saying that the function $z \mapsto {\bf f}_\varphi(z,g_f)$ is holomorphic for all $g_f \in \GSp_{2n}(\A_{\F,f})$ (cf.\,\cite[Lemma 7]{AS2001}). In particular, we have the Fourier expansion
\begin{align*}
{\bf f}_\varphi(z,g_f) &= \sum_{B \in {\rm Sym}_n(\F)}{\bf a}_B(g_f,{\bf f}_\varphi)e^{2\pi\sqrt{-1}\,\sum_{v \in S_\infty}{\rm tr}(\iota_v(B)z_v)},\\
{\bf a}_B(g_f,{\bf f}_\varphi) & = e^{2\pi {\rm tr}_{\F/\Q} ({\rm tr}(B))} \cdot W_{B} (g_f,\varphi).
\end{align*}
Therefore, under the assumption in Theorem \ref{T:q-expansion}, for all $\sigma \in {\rm Aut}(\C)$, ${}^\sigma\!\varphi$ is the unique holomorphic automorphic form of weight $\ell$ and has central character ${}^\sigma\!\chi$ such that
\[
\sigma\left( {\bf a}_B(g_f,{\bf f}_\varphi)\right) = {\bf a}_B(g_f,{\bf f}_{{}^\sigma\!\varphi})
\]
for all $B \in {\rm Sym}_n(\F)$ and $g_f \in \GSp_{2n}(\A_{\F,f})$.

\subsection{Degenerate principal series representations}\label{SS:degenerate p.s}

Let $\chi$ and $\mu$ be Hecke characters of $\A_\F^\times$.
We assume that $\chi$ is unitary.
For $s \in \C$, let 
\[
I(\chi \rtimes \mu,s) = {\rm Ind}_{P_n(\A_\F)}^{\GSp_{2n}(\A_\F)}\left(\chi|\mbox{ }|_{\A_\F}^{s-(n+1)/2}\rtimes \mu|\mbox{ }|_{\A_\F}^{n(n+1)/4}\right)
\]
be the degenerate principal series representation consisting of smooth and right $({\rm U}(n)^{\Sigma_\F}\times \GSp_{2n}(\hat{\o}_\F))$-finite functions $f:\GSp_{2n}(\A_\F)\rightarrow \C$ such that
\[
f({\bf m}(a,\nu){\bf n}(x) g) = \mu(\nu)\chi(\det a) |\det a|_{\A_\F}^{s}\cdot f(g).
\]
for $\nu \in \A_\F^\times$, $a \in \GL_n(\A_\F)$, $x \in {\rm Sym}_n(\A_\F)$, and $g \in \GSp_{2n}(\A_\F)$.
A function
\[
\C \times \GSp_{2n}(\A_\F) \longrightarrow \C,\quad (s,g)\longmapsto f^{(s)}(g)
\]
is called a holomorphic section of $I(\chi\rtimes\mu,s)$ if it satisfies the following conditions:
\begin{itemize}
\item For each $s \in \C$, the function $g \mapsto f^{(s)}(g)$ belongs to $I(\chi\rtimes\mu,s)$.
\item For each $g \in \GSp_{2n}(\A_\F)$, the function $s \mapsto f^{(s)}(g)$ is holomorphic.
\item $f^{(s)}$ is right $({\rm U}(n)^{\Sigma_\F}\times \GSp_{2n}(\hat{\o}_\F))$-finite.
\end{itemize}
A function $f^{(s)}$ on $\C \times \GSp_{2n}(\A_\F)$ is called a meromorphic section of $I(\chi\rtimes\mu,s)$ if there exists a non-zero entire function $\beta$ such that $\beta(s)f^{(s)}$ is a holomorphic section.

Let $v$ be a place of $\F$. Let $I(\chi_v\rtimes\mu_v,s)$ be the degenerate principal series representation of $\GSp_{2n}(\F_v)$ defined similarly as above. Holomorphic sections and meromorphic sections of $I(\chi_v\rtimes\mu_v ,s)$ are defined in similar way. 
Assume $v$ is finite. 
A meromorphic section $f_v^{(s)}$ is called a rational section if for any $g \in \GSp_6(\F_v)$, the function $s \mapsto f_v^{(s)}(g)$ is a rational function in $q_v^{-s}$. For a rational section $f_v^{(s)}$ and $\sigma \in {\rm Aut}(\C)$, we define the rational section ${}^\sigma\!f_v^{(s)}$ of $I({}^\sigma\!\chi_v\rtimes{}^\sigma\!\mu_v,s)$ by
\begin{align}\label{E:Galois action section}
{}^\sigma\!f_v^{(s)}\left({\bf m}(a,\nu){\bf n}(x) k\right) = {}^\sigma\!\mu_v(\nu){}^\sigma\!\chi_v(\det a) |\det a|_v^{s}\cdot {}^\sigma\!(f_v^{(s)}(k))
\end{align}
for $\nu \in \F_v^\times$, $a \in \GL_n(\F_v)$, $x\in{\rm Sym}_n(\F_v)$, and $k \in \GSp_{2n}(\o_{\F_v})$.
Here ${}^\sigma\!\omega(x) = \sigma (\omega(x))$ for a character $\omega$ of $\F_v^\times$. 

\subsection{Siegel Eisenstein series}

For a holomorphic section $f^{(s)}$ of $I(\chi \rtimes \mu,s)$, we define the Siegel Eisenstein series $E(f^{(s)})$ on $\GSp_{2n}(\A_\F)$ by the absolutely convergent series
\[
E(g,f^{(s)}) = \sum_{\gamma \in P_n(\F)\backslash \GSp_{2n}(\F)}f^{(s)}(\gamma g)
\]
for ${\rm Re}(s)>n+1$, and by meromorphic continuation otherwise.
We have the Fourier expansion
\[
E(g,f^{(s)}) = \sum_{B \in {\rm Sym}_n(\F)}W_B(g,E(f^{(s)})).
\]
For brevity, we write $W_B(f^{(s)}) = W_B(E(f^{(s)}))$ for the $B$-the Fourier coefficient.

Let $v$ be a place of $\F$. Let $\psi_v$ be the $v$-component of $\psi_\F$.
Let $f_v^{(s)}$ be a meromorphic section of $I(\chi_v \rtimes \mu_v,s)$. For $B_v \in {\rm Sym}_n(\F_v)$ with $\det B_v \neq 0$, we defined the degenerate Whittaker functional $W_{B_v}(f_v^{(s)}): \GSp_{2n}(\F_v) \rightarrow \C$ by
\[
W_{B_v}(g,f_v^{(s)}) = \int_{{\rm Sym}_n(\F_v)} f_v^{(s)}(J_n{\bf n}(x_v) g) \overline{\psi_v({\rm tr}(B_vx_v))}\,dx_v.
\]
Here $dx_v$ is the product measure on $\F_v^{n(n+1)/2}$ with respect to the isomorphism 
\[
{\rm Sym}_n(\F_v) \longrightarrow \F_v^{n(n+1)/2},\quad x_v \longmapsto (x_{v,ij})_{1\leq i \leq j \leq n};
\]
and the Haar measure on $\F_v$ is normalized so that ${\rm vol}(\o_{\F_v})=1$ if $v$ is finite, and we take the Lebesgue measure on $\F_v$ if $v \in S_\infty$.
When $f_v^{(s)}$ is holomorphic, the integral converges absolutely for ${\rm Re}(s)>n$ and admits holomorphic continuation to $\C$.
Moreover, when $v$ is finite and $f_v^{(s)}$ is holomorphic and rational, $W_{B_v}(g,f_v^{(s)})$ defines a polynomial in $\C[q_v^{s},q_v^{-s}]$ by the result of Karel \cite[Corollary 3.6.1]{Karel1979}. 
In particular, for a rational section $f_v^{(s)}$, $W_{B_v}(g,f_v^{(s)})$ defines a rational function in $q_v^{-s}$ for each $g \in \GSp_{2n}(\F_v)$.

\begin{lemma}\label{L:spherical Whittaker}
Assume $v$ is a finite place, $v \nmid 2$, $\chi_v,\mu_v,\psi_v$ are unramified, and $B_v \in {\rm Sym}_n(\o_{\F_v})$ with $\det B_v \in \o_{\F_v}^\times$. Let $f_{v,\circ}^{(s)}$ be the $\GSp_{2n}(\o_{\F_v})$-invariant good section of $I(\chi_v \rtimes \mu_v,s)$ normalized so that
\[
f_{v,\circ}^{(s)}(1) = L(s,\chi_v)\prod_{j=1}^{\lfloor \tfrac{n}{2} \rfloor}L(2s-2j,\chi_v^2).
\]
We have
\begin{align*}
W_{B_v}(1,f_{v,\circ}^{(s)}) = 
\begin{cases}
L(s-\tfrac{n}{2},\chi_v\chi_{B_v}) & \mbox{ if $n$ is even},\\
1 & \mbox{ if $n$ is odd}.
\end{cases}
\end{align*}
Here $\chi_{B_v}$ is the quadratic character of $\F_v^\times$ associated to $\F_v(\sqrt{(-1)^{n/2}\det B_v})/\F_v$ if $n$ is even.
\end{lemma}

\begin{proof}
We refer to \cite[Theorem 13.6]{Shimura1997} for the computation.
\end{proof}

\begin{lemma}\label{L:Galois conjugate Whittaker}
Assume $v$ is a finite place lying over a rational prime $p$. Let $f_v^{(s)}$ be a rational section of $I(\chi_v \rtimes \mu_v,s)$.
We have
\[
{}^\sigma\!\left(W_{B_v}\left( \bp u_{\sigma,p}^{-1}{\bf 1}_n & 0 \\ 0 & {\bf 1}_n\ep g,f_v^{(s)}\right)\right) = {}^\sigma\!\mu_v(u_{\sigma,p})^{-1}\cdot W_{B_v}(g,{}^\sigma\!f_v^{(s)})
\]
as rational functions in $q_v^{-s}$ for all $\sigma \in {\rm Aut}(\C)$ and $g \in \GSp_{2n}(\F_v)$.
Here $u_{\sigma,p} \in \Z_p^\times$ is the unique element such that $\sigma(\psi_v(x)) = \psi_v(u_{\sigma,p}x)$ for $x\in\F_v$.
\end{lemma}

\begin{proof}
By the result of Karel \cite[Theorem 3.6]{Karel1979}, the integral defining $W_{B_v}(f_v^{(s)})$ is actually a stable integral, that is, for all sufficiently large integer $N$ we have
\[
W_{B_v}(g,f_v^{(s)}) = \int_{{\rm Sym}_n(\varpi_v^{-N}\o_{\F_v})}f_v^{(s)}(J_n{\bf n}(x_v) g) \overline{\psi_v({\rm tr}(B_vx_v))}\,dx_v.
\]
Here $\varpi_v$ is a uniformizer of $\F_v$.
Let $\sigma \in {\rm Aut}(\C)$. Since ${\rm Sym}_n(\varpi_v^{-N}\o_{\F_v})$ is compact, for all sufficiently large integer $N$ we have
\begin{align*}
&{}^\sigma\!\left(W_{B_v}(g,f_v^{(s)})\right)\\
& = \int_{{\rm Sym}_n(\varpi_v^{-N}\o_{\F_v})}{}^\sigma\!f_v^{(s)}(J_n{\bf n}(x_v)g) \sigma\left(\overline{\psi_v({\rm tr}(B_vx_v))}\right)\,dx_v\\
& = \int_{{\rm Sym}_n(\varpi_v^{-N}\o_{\F_v})}{}^\sigma\!f_v^{(s)}(J_n{\bf n}(x_v) g) \overline{\psi_v({\rm tr}(u_{\sigma,p}B_vx_v))}\,dx_v\\
& = \int_{{\rm Sym}_n(\varpi_v^{-N}\o_{\F_v})}{}^\sigma\!f_v^{(s)}\left(\bp {\bf 1}_n & 0 \\ 0 & u_{\sigma,p}^{-1}{\bf 1}_n\ep J_n{\bf n}(x_v) \bp u_{\sigma,p}{\bf 1}_n & 0 \\ 0 & {\bf 1}_n\ep g\right) \overline{\psi_v({\rm tr}(B_vx_v))}\,dx_v\\
& = {}^\sigma\!\mu_v(u_{\sigma,p})^{-1}\cdot W_{B_v}\left(\bp u_{\sigma,p}{\bf 1}_n & 0 \\ 0 & {\bf 1}_n\ep g,{}^\sigma\!f_v^{(s)}\right)
\end{align*}
as rational functions in $q_v^{-s}$.
This completes the proof.
\end{proof}

\begin{lemma}\label{L:archimedean Whittaker}
Assume $v \in S_\infty$ and $\chi_v$ is trivial on $\R_{>0}^\times$. For a positive integer $\ell$ such that $(-1)^{\ell} = \chi_v(-1)$, let $f_{v,\ell}^{(s)}$ be the holomorphic section of $I(\chi_v \rtimes \mu_v,s)$ defined by
\begin{align}\label{E:archimedean section}
f_{v,\ell}^{(s)}(g) = {\rm sgn}^{n\ell}(\nu(g))\mu_v(\nu(g))|\nu(g)|_v^{ns}\cdot\det(c\,\sqrt{-1}+d)^{-\ell}\left\vert \det(c\,\sqrt{-1}+d)\right\vert^{\ell-s}
\end{align}
for $g = \bp a & b \\ c & d\ep \in \GSp_{2n}(\F_v)$.
Then for $a_v \in \GL_n(\F_v)$, we have
\begin{align*}
&W_{B_v}({\bf m}(a_v),f_{v,\ell}^{(s)})\vert_{s=\ell}\\
&=
(4\pi)^{-n(n-1)/4}\prod_{j=0}^{n-1}\Gamma(\ell-\tfrac{j}{2})^{-1}(2\pi\sqrt{-1})^{n\ell}\cdot (\det B_v)^{\ell-(n+1)/2}(\det a_v)^{\ell}e^{-2\pi{\rm tr}({}^t\!a_vB_va_v)}
\end{align*}
if $B_v$ is positive definite, and equal to zero otherwise.
\end{lemma}

\begin{proof}
The assertion is a direct consequence of the computation of Shimura in \cite[(4.34K) and (4.35K)]{Shimura1982}.
\end{proof}

Assume now that $\chi$ and $\mu$ are algebraic Hecke characters and $\chi$ has parallel signature ${\rm sgn}(\chi)$. 
For $\ell\geq1$ such that $(-1)^\ell = {\rm sgn}(\chi)$, let $f_{\infty,\ell}^{(s)}$ be the holomorphic section of $\bigotimes_{v \in S_\infty}I(\chi_v\rtimes\mu_v,s)$ defined by
\[
f_{\infty,\ell}^{(s)} = \bigotimes_{v \in S_\infty} f_{v,\ell}^{(s)}.
\]
The algebraicity of holomorphic Eisenstein series in the region of absolute convergence was proved in \cite[Theorem 3.4.7 and Appendix A.2]{GH1993} for $n=3$ based on the result of Harris \cite{Harris1984}. The result for general $n$ can be proved in a similar way.

\begin{thm}[Harris]\label{T:Harris}
Let $f^{(s)} = \bigotimes_{v \nmid \infty} f_v^{(s)}$ be a meromorphic section of $\bigotimes_{v \nmid \infty}I(\chi_v\rtimes\mu_v,s)$ and $\ell > n+1$ with $(-1)^\ell = {\rm sgn}(\chi)$.
Assume the following conditions are satisfied:
\begin{itemize}
\item $f_v^{(s)}$ is a rational section and holomorphic for ${\rm Re}(s)>n+1$ for all $v$.
\item $f_v^{(s)} = f_{v,\circ}^{(s)}$ for almost all $v$.
\end{itemize}
Then the following assertions hold:
\begin{itemize}
\item[(1)] The Eisenstein series $E(f_{\infty,\ell}^{(s)}\otimes f^{(s)})$ is holomorphic at $s=\ell$.
\item[(2)] The automorphic form $E^{[\ell]}(f^{(s)}) = E(f_{\infty,\ell}^{(s)}\otimes f^{(s)})\vert_{s=\ell}$ is holomorphic of weight $\ell$ and has central character $\chi^n\mu^2|\mbox{ }|_{\A_\F}^{n\ell}$. 
\item[(3)] For $\sigma \in {\rm Aut}(\C)$, we have
\begin{align}\label{E:Galois equiv. E.S. even}
\begin{split}
&\frac{{}^\sigma\!E^{[\ell]}(f^{(s)})}{\sigma\left(|D_\F|^{(n+2)/4}\cdot(2\pi\sqrt{-1})^{(n\ell+\ell-n(n+2)/4)d}\cdot G(\chi\cdot\mu^{-1})\right)}\\
& = \frac{E^{[\ell]}({}^\sigma\!f^{(s)})}{|D_\F|^{(n+2)/4}\cdot(2\pi\sqrt{-1})^{(n\ell+\ell-n(n+2)/4)d}\cdot G({}^\sigma\!\chi\cdot{}^\sigma\!\mu^{-1})}
\end{split}
\end{align}
if $n$ is even; and
\begin{align}\label{E:Galois equiv. E.S.}
\begin{split}
&\frac{{}^\sigma\!E^{[\ell]}(f^{(s)})}{\sigma\left(|D_\F|^{(n+1)/4}\cdot(2\pi\sqrt{-1})^{(n\ell-(n^2-1)/4)d}\cdot G(\mu^{-1})\right)}\\
& = \frac{E^{[\ell]}({}^\sigma\!f^{(s)})}{|D_\F|^{(n+1)/4}\cdot(2\pi\sqrt{-1})^{(n\ell-(n^2-1)/4)d}\cdot G({}^\sigma\!\mu^{-1})}
\end{split}
\end{align}
if $n$ is odd.
\end{itemize}
\end{thm}

In the following proposition, we extend the algebraicity to $\tfrac{n}{2}<\ell \leq n+1$  under certain local assumption.
For a finite place $v$ and a Schwartz function $\Phi_v \in \mathcal{S}({\rm Sym}_3(\F_v))$, let $f_{\itPhi_v}^{(s)}$ be the rational section of $I(\chi_v\rtimes\mu_v,s)$ such that $f_{\itPhi_v}^{(s)}$ is supported in $P_3(\F_v)J_3P_3(\F_v)$ and
\[
f_{\itPhi_v}^{(s)}(J_3{\bf n}(x)) = \itPhi_v(x).
\]
We are inspired by the construction of Eisenstein series due to Shimura \cite[\S\,18]{Shimura1997}, where the local sections at bad places are of the form $f_{\Phi}^{(s)}$.

\begin{prop}\label{P:holomorphy}
Let $f^{(s)} = \bigotimes_{v \nmid \infty} f_v^{(s)}$ be a meromorphic section of $\bigotimes_{v \nmid \infty}I(\chi_v\rtimes\mu_v,s)$ and $\ell > \tfrac{n}{2}$ with $(-1)^\ell = {\rm sgn}(\chi)$.
Assume the following conditions are satisfied:
\begin{itemize}
\item When $\F=\Q$, $n$ is even, and $\ell=1+\tfrac{n}{2}$, we have $\chi^2 \neq 1$.
\item $f_v^{(s)}$ is a rational section and holomorphic for ${\rm Re}(s)>\tfrac{n}{2}$ for all $v$.
\item $f_v^{(s)} = f_{v,\circ}^{(s)}$ for almost all $v$.
\item $f_{v_0}^{(s)} = f_{\itPhi_{v_0}}^{(s)}$ for some $v_0$ and $\itPhi_{v_0} \in \mathcal{S}({\rm Sym}_n(\F_{v_0}))$, and the Fourier transform 
\[
\widehat{\itPhi}_{v_0}(x_{v_0}) = \int_{{\rm Sym}_3(\F_{v_0})}\Phi_{v_0}(y_{v_0})\psi_{v_0}({\rm tr}(x_{v_0}y_{v_0}))\,dy_{v_0}
\]
is supported in ${\rm Sym}_n(\F_{v_0}) \cap \GL_n(\F_{v_0})$.
\end{itemize}
Then the assertions (1)-(3) in Theorem \ref{T:Harris} hold.
\end{prop}

\begin{proof}
Fix an open compact subgroup $K_{v_0}$ of $\GSp_{2n}(\F_{v_0})$ such that $f_{v_0}^{(s)}$ is right invariant by $K_{v_0}$. 
By the Bruhat decomposition for $P_n\backslash \GSp_{2n}/P_n$ and the condition on the support of $f_{v_0}^{(s)}$, for ${\rm Re}(s)>n+1$ and $B \in {\rm Sym}_n(\F)$, the $B$-th Fourier coefficient $W_B(g,f_{\infty,\ell}^{(s)}\otimes f^{(s)})$ of the Eisenstein series $E(g,f_{\infty,\ell}^{(s)}\otimes f^{(s)})$ is equal to
\begin{align*}
W_B(g,f_{\infty,\ell}^{(s)}\otimes f^{(s)}) = \int_{{\rm Sym}_n(\A_\F)} f^{(s)}(J_n{\bf n}(x) g) \overline{\psi({\rm tr}(Bx))}\,dx^{\rm Tam}.
\end{align*}
The second condition on $f_{v_0}^{(s)}$ implies that if $\det B=0$, then
\[
\int_{{\rm Sym}_n(\F_{v_0})} f_{v_0}^{(s)}(J_n{\bf n}(x_{v_0}) g_{v_0}) \overline{\psi_{v_0}({\rm tr}(Bx_{v_0}))}\,dx_{v_0} =0
\]
for all $g_{v_0} \in P_n(\F_{v_0}) K_{v_0}$.
In particular, for $g \in \GSp_{2n}(\A_\F)$ such that $g_{v_0} \in P_n(\F_{v_0})K_{v_0}$ and $B \in {\rm Sym}_n(\F)$ with $\det B =0$, we have $W_B(g,f_{\infty,\ell}^{(s)}\otimes f^{(s)})=0$ for ${\rm Re}(s)>n+1$. Therefore, we have the Fourier expansion
\begin{align}\label{E:holomorphy proof 1}
E(g,f_{\infty,\ell}^{(s)}\otimes f^{(s)}) = \sum_{B \in {\rm Sym}_n(\F),\,\det B \neq 0}W_B(g,f_{\infty,\ell}^{(s)}\otimes f^{(s)})
\end{align}
for all $g \in \GSp_{2n}(\A_\F)$ such that $g_{v_0} \in P_n(\F_{v_0})K_{v_0}$ and ${\rm Re}(s)>n+1$. 
Let
\[
\frak{S} = \GSp_{2n}(\A_\F)\bigcap \left(\GSp_{2n}^+(\R)^{S_\infty}\cdot \prod_{v \nmid \infty,\, v\neq v_0}\GSp_{2n}(\F_v)\cdot P_n(\F_{v_0})K_{v_0}\right).
\]
We compute the non-degenerate Fourier coefficients of $E(g,f_{\infty,\ell}^{(s)}\otimes f^{(s)})$ for $g \in \frak{S}$.
Since $\GSp_{2n}^+(\R) = \R^\times P_n(\R) {\rm U}(n)$, it suffices to consider $g \in \frak{S}$ with $g_\infty = {\bf m}(a_\infty) \in P_n(\F_\infty)$ for some $a_\infty \in \GL_n(\R)^{S_\infty}$.
By Lemmas \ref{L:spherical Whittaker} and \ref{L:archimedean Whittaker}, for $B \in {\rm Sym}_{2n}(\F)$ with $\det B \neq 0$, we have
\begin{align}\label{E:Fourier coefficient}
\begin{split}
& W_B(g,f_{\infty,\ell}^{(s)}\otimes f^{(s)})\vert_{s=\ell}\\
& = (4\pi)^{-n(n-1)d/4}\prod_{j=0}^{n-1}\Gamma(\ell-\tfrac{j}{2})^{-d}(2\pi\sqrt{-1})^{n\ell d}\cdot {\rm N}_{\F/\Q}(\det B)^{\ell-(n+1)/2}(\det a_\infty)^{\ell}\cdot e^{2\pi\sum_{v \in S_\infty}{\rm tr}({}^t\!a_v\iota_v(B)a_v)}\\
&\times |D_\F|^{-n(n+1)/4}\prod_{v \in S \setminus S_\infty}W_B(g_v,f_v^{(s)})\vert_{s=\ell}\cdot \begin{cases}L^S(\ell-\tfrac{n}{2},\chi\chi_B) & \mbox{ if $n$ is even},\\
1 & \mbox{ if $n$ is odd}.
\end{cases}
\end{split}
\end{align}
Here $S$ is a sufficiently large finite set of places containing $S_\infty$.
Note that the condition $(-1)^\ell = {\rm sgn}(\chi)$ implies that $L^S(s,\chi\chi_B)$ is holomorphic at $s = \ell-\tfrac{n}{2}$.
Also the factor $|D_\F|^{-n(n+1)/4}$ is due to the comparison between the Tamagawa measure $dx^{\rm Tam}$ and the product measure $\prod_v dx_v$ on ${\rm Sym}_n(\A_\F)$.
On the other hand, when $B$ is not totally positive definite, we have $W_B(g_f,f_{\infty,\ell}^{(s)}\otimes f^{(s)})\vert_{s=\ell} =0$.
Indeed, if either $\ell>1+\tfrac{n}{2}$ or $\chi\neq\chi_B$, then $\prod_{v \nmid \infty}W_B(g_v,f_v^{(s)})$ is holomorphic at $s = \ell$.
Thus $W_B(g_f,f_{\infty,\ell}^{(s)}\otimes f^{(s)})\vert_{s=\ell}$ vanishes as $\prod_{v \in S_\infty}W_B({\bf m}(a_v),f_{v,\ell}^{(s)})$ has a zero at $s = \ell$.
If $\ell = 1+\tfrac{n}{2}$ and $\chi=\chi_B$, then $\prod_{v \nmid \infty}W_B(g_v,f_v^{(s)})$ has at most a simple pole at $s=\ell$.
In this case, $\F\neq\Q$ by the first condition, hence $\prod_{v \in S_\infty}W_B({\bf m}(a_v),f_{v,\ell}^{(s)})$ has a zero of order at least $2$ at $s=\ell$.
Thus $W_B(g_f,f_{\infty,\ell}^{(s)}\otimes f^{(s)})\vert_{s=\ell}$ vanishes.
In conclusion, the right-hand side of (\ref{E:holomorphy proof 1}) is holomorphic at $s=\ell$ for all $g \in \frak{S}$. In particular, $E(f_{\infty,\ell}^{(s)}\otimes f^{(s)})$ is holomorphic at $s  = \ell$.

For the second assertion, consider the function ${\bf E}(f^{(s)}) : \frak{H}^{S_\infty}\times\GSp_{2n}(\A_{\F,f}) \rightarrow \C$ defined as in (\ref{E:H-S modular form}) associated to $E^{[\ell]}(f^{(s)})$.
Then ${\bf E}(f^{(s)})$ is real analytic on $\frak{H}^{S_\infty}$.
It is clear from the definition of $f_{\infty,\ell}^{(s)}$ that 
\[
E^{[\ell]}(f^{(s)})(gak_\infty) = \chi^n\mu^2|\mbox{ }|_{\A_\F}^{n\ell}(a)(\det k_\infty)^\ell \cdot E^{[\ell]}(f^{(s)})(g)
\]
for all $a \in \A_\F^\times$, $k_\infty \in {\rm U}(n)^{S_\infty}$, and $g \in \GSp_{2n}(\A_\F)$.
It remains to show that $(\frak{p}^-)^{S_\infty}\cdot E^{[\ell]}(f^{(s)}) =0$. As we recalled in the last paragraph of \S\,\ref{S:H-S modular form}, it is equivalent to showing that the function $z \mapsto {\bf E}(f^{(s)})(z,g_f)$ is holomorphic for all $g_f \in \GSp_{2n}(\A_{\F,f})$.
We have the Fourier expansion
\begin{align*}
{\bf E}(f^{(s)})(z,g_f) &= \sum_{B \in {\rm Sym}_n(\F)}{\bf a}_B(y,g_f,{\bf E}(f^{(s)}))e^{2\pi\sqrt{-1}\,\sum_{v \in S_\infty}{\rm tr}(\iota_v(B)z_v)},\quad z=x+y\sqrt{-1} \in \frak{H}_n^{S_\infty}\\
{\bf a}_B(y,g_f,{\bf E}(f^{(s)})) & = e^{2\pi\sum_{v \in S_\infty}{\rm tr}(\iota_v(B)y_v)}(\det y)^{-\ell/2}\cdot W_B\left.\left({\bf m}(y^{1/2})\cdot g_f,f_{\infty,\ell}^{(s)}\otimes f^{(s)}\right)\right\vert_{s=\ell}.
\end{align*} 
By (\ref{E:Fourier coefficient}), for all $g_f \in \GSp_{2n}(\A_{\F,f}) \cap \frak{S}$, we have ${\bf a}_B(y,g_f,{\bf E}(f^{(s)})) ={\bf a}_B(1,g_f,{\bf E}(f^{(s)}))$ and ${\bf a}_B(1,g_f,{\bf E}(f^{(s)}))=0$ if $B$ is not totally positive definite. Therefore, $X\cdot E^{[\ell]}(f^{(s)})(g)=0$ for all $X \in (\frak{p}^-)^{S_\infty}$ and $g \in \frak{S}$.
We conclude from the strong approximation theorem for $\Sp_{2n}$ that $(\frak{p}^-)^{S_\infty}\cdot E^{[\ell]}(f^{(s)})=0$.

Now we prove the third assertion.
Let $B \in {\rm Sym}_n(\F)$ be totally positive definite. When $n$ is even, by \cite[Proposition 3.1]{Shimura1978}, we have
\[
\sigma \left( \frac{L^S(\ell-\tfrac{n}{2},{}^\sigma\!\chi\chi_B)}{|D_\F|^{1/2}\cdot (2\pi\sqrt{-1})^{(\ell-n/2)d}\cdot G({}^\sigma\!\chi\chi_B)}\right) = \frac{L^S(\ell-\tfrac{n}{2},\chi\chi_B)}{|D_\F|^{1/2}\cdot (2\pi\sqrt{-1})^{(\ell-n/2)d}\cdot G(\chi\chi_B)}
\]
for $\sigma \in {\rm Aut}(\C)$. Also note that
\[
\frac{G(\chi_B)}{{\rm N}_{\F/\Q}(\det B)\cdot (\sqrt{-1})^{nd/2}} \in \Q^\times
\]
if $n$ is even, and 
\[
(4\pi)^{n(n-1)/4}\prod_{j=0}^{n-1}\Gamma(\ell-\tfrac{j}{2}) \in \begin{cases}
\pi^{n^2/4}\cdot \Q^\times & \mbox{ if $n$ is even},\\
\pi^{(n^2-1)/4}\cdot \Q^\times & \mbox{ if $n$ is odd}.
\end{cases}
\]
Together with Lemma \ref{L:Galois conjugate Whittaker} for $v \in S \setminus S_\infty$ and property (\ref{E:Galois Gauss sum}) for Gauss sum, we deduce that 
\begin{align}\label{E:holomorphy proof 2}
\begin{split}
&\sigma\left( e^{2\pi {\rm tr}_{\F/\Q}({\rm tr}(B))}\cdot\frac{W_B\left(\left.\bp u^{-1}{\bf 1}_n & 0 \\ 0 & {\bf 1}_n\ep g_f,f_{\infty,\ell}^{(s)}\otimes f^{(s)}\right)\right\vert_{s = \ell}}{|D_\F|^{1/2+n(n+1)/4}\cdot(2\pi\sqrt{-1})^{(n\ell+\ell-n(n+2)/4)d}\cdot G(\chi\cdot\mu^{-1})}\right) \\
&= e^{2\pi {\rm tr}_{\F/\Q}({\rm tr}(B))}\cdot\frac{W_B(g_f,f_{\infty,\ell}^{(s)}\otimes {}^\sigma\!f^{(s)}) \vert_{s = \ell}}{|D_\F|^{1/2+n(n+1)/4}\cdot(2\pi\sqrt{-1})^{(n\ell+\ell-n(n+2)/4)d}\cdot G({}^\sigma\!\chi\cdot{}^\sigma\!\mu^{-1})}
\end{split}
\end{align}
if $n$ is even; and
\begin{align}\label{E:holomorphy proof 3}
\begin{split}
&\sigma\left( e^{2\pi {\rm tr}_{\F/\Q}({\rm tr}(B))}\cdot\frac{W_B\left(\left.\bp u^{-1}{\bf 1}_n & 0 \\ 0 & {\bf 1}_n\ep g_f,f_{\infty,\ell}^{(s)}\otimes f^{(s)}\right)\right\vert_{s = \ell}}{|D_\F|^{n(n+1)/4}\cdot(2\pi\sqrt{-1})^{(n\ell-(n^2-1)/4)d}\cdot G(\mu^{-1})}\right) \\
&= e^{2\pi {\rm tr}_{\F/\Q}{\rm tr}(B)}\cdot\frac{W_B(g_f,f_{\infty,\ell}^{(s)}\otimes {}^\sigma\!f^{(s)}) \vert_{s = \ell}}{|D_\F|^{n(n+1)/4}\cdot(2\pi\sqrt{-1})^{(n\ell-(n^2-1)/4)d}\cdot G({}^\sigma\!\mu^{-1})}
\end{split}
\end{align}
if $n$ is odd, for $\sigma \in {\rm Aut}(\C)$.
Here $u_\sigma \in \widehat{\Z}^\times$ is the unique element such that $\sigma(\psi_\Q(x)) = \psi_\Q(u_\sigma x)$ for $x \in \A_{\Q,f}$.
Comparing (\ref{E:holomorphy proof 2}) and (\ref{E:holomorphy proof 3}) with Theorem \ref{T:q-expansion}, we conclude that (\ref{E:Galois equiv. E.S. even}) and (\ref{E:Galois equiv. E.S.}) hold for all $g \in \frak{S}$. It then follows from the strong approximation theorem for $\Sp_{2n}$ that the (\ref{E:Galois equiv. E.S. even}) and (\ref{E:Galois equiv. E.S.}) hold for all $g \in \GSp_{2n}(\A_\F)$.
This completes the proof.
\end{proof}


\section{Triple product $L$-functions}\label{S:triple L-function}\label{S:triple}

Let $\E$ be a totally real \'etale cubic algebra over $\F$.
Let $(\itPi,V_\itPi)$ be an irreducible cuspidal automorphic representation of $({\rm R}_{\E/\F}\GL_2)(\A_\F) = \GL_2(\A_\E)$ with central character $\omega_\itPi$. Put $\omega = \omega_\itPi\vert_{\A_\F^\times}$. In this section, we recall the definition of the triple product $L$-function $L(s,\itPi,{\rm As})$ and Garrett's integral representation from adelic point of view following \cite{PSR1987}. 

\subsection{Local factors via the Weil--Deligne representations}

Let $v$ be a place of $\F$ and $\psi_v$ a non-trivial additive character of $\F_v$. Let $L_{\F_v}$ be the Weil--Deligne group of $\F_v$.
We identify the Langlands dual group ${}^L({\rm R}_{\E_v/\F_v}\GL_2)$ of ${\rm R}_{\E_v/\F_v}\GL_2$ with $\GL_2(\C)^3 \rtimes \Gal(\overline{F}_v/F_v)$ (cf.\,\cite[\S\,5]{Borel1979}), where the action of $\Gal(\overline{F}_v/F_v)$ on $\GL_2(\C)^3$ is described as follows:
We have the following three cases
\begin{align*}
\begin{cases}
\E_v=\F_v\times \F_v \times \F_v & \mbox{ Case 1},\\
\E_v=\F_v'\times \F_v \mbox{ for some quadratic extension $\F_v'$ of $\F_v$}& \mbox{ Case 2},\\
\E_v \mbox{ is a field} & \mbox{ Case 3}.
\end{cases}
\end{align*}
In Case 1, the action is trivial. In Case 2, the action is the permutation on the first two copies of $\GL_2(\C)^3$ through the natural surjection $\Gal(\overline{\F}_v/\F_v) \rightarrow \Gal(\F_v'/\F_v)$. In Case, the action is the permutation on $\GL_2(\C)^3$ through the natural surjection $\Gal(\overline{\F}_v/\F_v) \rightarrow \Gal(\widetilde{\E}_v/\F_v)$, where $\widetilde{\E}_v$ is the Galois closure of $\E_v/\F_v$. Let
\[
{\rm As} : {}^L({\rm R}_{\E_v/\F_v}\GL_2) \longrightarrow \GL(\C^2\otimes\C^2\otimes\C^2)
\]
be the Asai cube representation defined by
\[
{\rm As}(g_1,g_2,g_3)\cdot(v_1\otimes v_2 \otimes v_3) = (g_1\cdot v_1,g_2\cdot v_2,g_3\cdot v_3)
\]
and the action of $\Gal(\overline{\F}_v/\F_v)$ on $\C^2\otimes\C^2\otimes\C^2$ is similarly to the one on $\GL_2(\C)^3$. 
Let $\itPi_v$ be the local component of $\itPi_v$ at $v$ and $\phi_{\itPi_v} : L_{\F_v} \rightarrow {}^L({\rm R}_{\E_v/\F_v}\GL_2)$ the corresponding Langlands parameter via the local Langlands correspondence. 
Then we have a $8$-dimensional admissible representation
\[
{\rm As}\circ\phi_{\itPi_v}: L_{\F_v}\longrightarrow \GL(\C^2\otimes\C^2\otimes\C^2).
\]
We denote by 
\[
L(s,\itPi_v,{\rm As}),\quad \varepsilon(s,\itPi_v,{\rm As},\psi_v)
\]
the $L$-factor and $\varepsilon$-factor associated to ${\rm As}\circ\phi_{\itPi_v}$ defined as in \cite[\S\,3]{Tate1979}.
For example, if $\F_v=\R$, $\E_v=\R^3$, $\psi_v(x)=e^{2\pi\sqrt{-1}\,x}$, and 
\[
\itPi_v = \boxtimes_{i=1}^3D_{\kappa_{i,v}}
\]
for some $\kappa_{1,v},\kappa_{2,v},\kappa_{3,v}\geq 1$ satisfying the balanced condition (\ref{E:balanced}), then we have
\begin{align}\label{E:local factor archimedean}
\begin{split}
L(s,\itPi_v,{\rm As}) &= \Gamma_\C(s+\tfrac{\kappa_{1,v}+\kappa_{2,v}+\kappa_{3,v}-3}{2})\Gamma_\C(s+\tfrac{\kappa_{1,v}+\kappa_{2,v}-\kappa_{3,v}-1}{2})\\
&\times\Gamma_\C(s+\tfrac{\kappa_{1,v}+\kappa_{3,v}-\kappa_{2,v}-1}{2})\Gamma_\C(s+\tfrac{\kappa_{2,v}+\kappa_{3,v}+\kappa_{1,v}-1}{2}),\\
\varepsilon(s,\itPi_v,{\rm As},\psi_v) &= (\sqrt{-1})^{2\kappa_{1,v}+2\kappa_{2,v}+2\kappa_{3,v}-3}.
\end{split}
\end{align}
Here $\Gamma_\C(s) = 2(2\pi)^{-s}\Gamma(s)$.
Also, if $v$ is finite, then for $\sigma \in {\rm Aut}(\C)$ we have
\begin{align}\label{E:local factor non-archimedean}
{}^\sigma\!L(s+\tfrac{1}{2},\itPi_v,{\rm As}) = L(s+\tfrac{1}{2},{}^\sigma\!\itPi_v,{\rm As}),\quad {}^\sigma\!\varepsilon(s+\tfrac{1}{2},\itPi_v,{\rm As},\psi_{v}) = {}^\sigma\!\omega_v(u_{\sigma,v})^4\cdot\varepsilon(s+\tfrac{1}{2},{}^\sigma\!\itPi_v,{\rm As},\psi_{v})
\end{align}
as rational functions in $q_v^{-s}$, where $u_{\sigma,v} \in \o_{\F_v}^\times$ is the unique element such that $\sigma(\psi_v(x)) = \psi_v(u_{\sigma,v} x)$ for $x \in \F_v$.
The above equalities can be easily verified using the results of Clozel \cite[Lemme 4.6]{Clozel1990} and Henniart \cite[Propri\'et\'e 3, \S\,7]{Henniart2001}.

\subsection{Garrett's integral representation}

We define the triple product $L$-function by the Euler product
\[
L(s,\itPi, {\rm As}) = \prod_{v}L(s,\itPi_v,{\rm As}).
\]
The Euler product converges absolutely when ${\rm Re}(s)$ sufficiently large and admits meromorphic continuation to $s \in \C$ by the results of Ikeda \cite{Ikeda1989}, \cite{Ikeda1992} and Piatetski-Shapiro--Rallis \cite{PSR1987}. They also establish functional equation for the triple product $L$-function. Combining with the results of Ramakrishinan \cite{Rama2000} and the author \cite{Chen2021b} on the comparison between gamma factors, we can now state the functional equation as follows:
Define the $\varepsilon$-factor $\varepsilon(s,\itPi,{\rm As})$ by
\[
\varepsilon(s,\itPi,{\rm As})=\prod_v\varepsilon(s,\itPi_v,{\rm As},\psi_v).
\]
Here $\psi=\bigotimes_v\psi_v$ is any non-trivial additive character of $\F\backslash\A_\F$. Then we have
\begin{align}\label{E:fe}
L(s,\itPi,{\rm As}) = \varepsilon(s,\itPi,{\rm As})L(1-s,\itPi^\vee,{\rm As}).
\end{align}
Here $\itPi^\vee$ is the contragredient of $\itPi$.

\subsubsection{Preliminaries}

Let $\G$ be the linear algebraic group over $\F$ defined by
\[
\G = \{g \in {\rm R}_{\E/\F}\GL_{2} \, \vert\, \det g \in \GL_1\}.
\]
Let $(V,\<\, ,\,\>)$ be the nondegenerate symplectic form over $\F$ defined by  
\[
V=({\rm R}_{\E/\F}\,{\mathbb G}_a)^2,\quad\<x,y\> = {\rm tr}_{\E/\F}(x_1y_2-x_2y_1)
\]
for $x = (x_1,x_2) , y= (y_1,y_2) \in V$. 
Let
\[
\GSp(V) = \{g \in {\rm R}_{\E/\F}\GL_{2} \, \vert\, \<xg,yg\> = \nu(g)\<x,y\> \mbox{ for all }x,y \in V,\, \nu(g) \in \GL_1\}
\]
be the similitude symplectic group associated to $\<\, ,\,\>$.
Note that ${\rm R}_{\E/\F}\GL_{2}$ acts on $V$ from the right.
It is easy to verify that $\G$ is a subgroup of $\GSp(V)$ and $\det g = \nu(g)$ for $g \in \G$.
Let $X$, $Y$, and $Y^0$ be maximal isotropic subspaces of $V$ defined by
\begin{align*}
X &= \{(x,0) \in V \,\vert\, x \in {\rm R}_{\E/\F}\,{\mathbb G}_a\},\quad Y = \{(0,y) \in V \,\vert\, y \in {\rm R}_{\E/\F}\,{\mathbb G}_a\},\\
\quad Y^0 &= \{ (x,y) \in V \, \vert\, x \in \mathbb{G}_a,\, {\rm tr}_{\E/\F}(y)=0 \}.
\end{align*}
Define an isomorphism between $Y(F)$ and $Y_0(F)$ by
\begin{align}\label{E:iso}
\begin{split}
Y(\F) \longrightarrow Y^0(\F),\quad (0,x)\longmapsto (3x,0),\quad
(0,y)\longmapsto (0,y)
\end{split}
\end{align}
for $x \in \F$ and $y \in \E$ with ${\rm tr}_{\E/\F}(y)=0$. 
We denote by $P$ and $R^0$ the stabilizers of $Y^0$ in $\GSp(V)$ and $\G$, respectively. Let $U^0$ be the unipotent radical of $R^0$.
Note that
\begin{align*}
R^0 &= \{ {\bf t}(a_1){\bf d}(a_2){\bf n}(x) \, \vert \, a_1,a_2 \in \GL_1,\, x \in {\rm R}_{\E/\F}\,\mathbb{G}_a,\, {\rm tr}_{\E/\F}(x)=0\},\\
U^0 &=  \{ {\bf n}(x) \, \vert \, x \in {\rm R}_{\E/\F}\,\mathbb{G}_a ,\, {\rm tr}_{\E/\F}(x)=0\}. 
\end{align*}
For $\nu \in \GL_1$ and $a \in \GL(X)$, let ${\bf m}(a,\nu) \in P$ defined by
\[
{\bf m}(a,\nu)\vert_{X} = a,\quad {\bf m}(a,\nu)\vert_{Y} = \nu\cdot (a^*)^{-1},
\]
where $a^* \in \GL(Y)$ is the unique linear transformation such that
\[
\<(x\cdot a,0),(0,y)\> = \<(x,0),(0,y\cdot a^*)\>.
\]
In the matrix form with respect to the complete polarization $V=X\oplus Y$, we have
\[
{\bf m}(a,\nu) = \bp a & 0 \\ 0 & \nu(a^*)^{-1}\ep.
\]
Then $M_P = \{{\bf m}(a,\nu) \,\vert\, a\in \GL(X),\,\nu\in\GL_1\}$ is a Levi component of $P$.
Fix an integral basis $\{x_1,x_2,x_3\}$ of $\E$ over $\F$ and $\{x_1^*,x_2^*,x_3^*\}$ be its dual basis. 
Put $e_i = (x_i,0)$ and $e_i^* = (0,x_i^*)$ for $i=1,2,3$.
Then $\{e_1,e_2,e_3,e_1^*,e_2^*,e_3^*\}$ is a symplectic basis of $V(\F)$ over $\F$, that is, we have
\[
\<e_i,e_j^*\> = \delta_{ij}.
\]
Let $K$ be the maximal open compact subgroup of $\GSp(V)(\A_\F)$ stabilizing the $\hat{\o}_{\F}$-lattice 
\[
\hat{\o}_{\F} e_1 \oplus \hat{\o}_{\F} e_2 \oplus \hat{\o}_{\F} e_3 \oplus \hat{\o}_{\F} e_1^* \oplus \hat{\o}_{\F} e_2^* \oplus \hat{\o}_{\F} e_3^* \subset V(\A_\F).
\]
For finite place $v$, let $K_v$ be the local component of $K$ at $v$.
For $v \in S_\infty$, let $w_{1,v},w_{2,v},w_{3,v}$ be the places of $\E$ lying over $v$. We identify $\E_{w_{i,v}}$ with $\R$ via the real embedding of $\E$ associated to $v$. 
We then identify $G(\F_v)$ with $(\GL_2(\R)^3)_0 = \{(g_1,g_2,g_3) \in \GL_2(\R)^3\,\vert\,\det g_1 = \det g_2 = \det g_3\}$ and $\GSp(V)(\F_v)$ with $\GSp_6(\R)$ in a natural way. In particular, the diagram
\[
\begin{tikzcd}
G(\F_v)   \arrow[r] \arrow[d] & \GSp(V)(\F_v)\arrow[d]\\
(\GL_2(\R)^3)_0 \arrow[r, "\iota"]  & \GSp_6(\R)
\end{tikzcd}
\]
commutes, where
\begin{align}\label{E:embedding}
\iota\left(\bp a_1 & b_1 \\ c_1 & d_1\ep,\bp a_2 & b_2 \\ c_2 & d_3\ep,\bp a_3 & b_3 \\ c_3 & d_3\ep\right) = 
\begin{pmatrix}
  a_1 & 0   & 0   & b_1 & 0   & 0   \\
  0   & a_2 & 0   & 0   & b_2 & 0   \\
  0   & 0   & a_3 & 0   & 0   & b_3 \\
  c_1 & 0   & 0   & d_1 & 0   & 0   \\
  0   & c_2 & 0   & 0   & d_2 & 0   \\
  0   & 0   & c_3 & 0   & 0   & d_3 \\ 
 \end{pmatrix}.
\end{align}
For $s \in \C$, let $I(\omega,s)$ be the degenerate principal series representation consisting of smooth and right $({\rm U}(3)^{S_\infty}\times K)$-finite function $f:\GSp(V)(\A_\F)\rightarrow \C$ such that $f$ is left-invariant by the unipotent radical of $P(\A_\F)$ and 
\[
f({\bf m}(a,\nu)g) = \omega^{-2}|\mbox{ }|_{\A_\F}^{-3s-3}(\nu)\cdot\omega|\mbox{ }|_{\A_\F}^{2s+2}(\det a)\cdot f(g)
\]
for ${\bf m}(a,\nu) \in M_P(\A_\F)$.
In the notation of \S\,\ref{SS:degenerate p.s}, we have
\begin{align}\label{E:degenerate p.s}
I(\omega,s) = I(\omega|\mbox{ }|_{\A_\F}^{-{\sf w}} \rtimes \omega^{-2}|\mbox{ }|_{\A_\F}^{-3s-3},2s+{\sf w}+2).
\end{align}

Let $dg=\prod_vdg_v$ be the Haar measure on $G(\A_\F)$ defined as follows:
For finite place $v$, $dg_v$ is normalized so that ${\rm vol}(\GL_2(\o_{\E_v})\cap G(\F_v),dg_v)=1$.  
For $v \in S_\infty$, we have
\[
dg_v = |a|_v^{-1}|a_1a_2a_3|_v^{-2}d a\prod_{i=1}^3 dx_i\,d a_i \,dk_i
\]
for $g_v = ({\bf t}(a){\bf n}(x_1){\bf m}(a_1)k_1,{\bf t}(a){\bf n}(x_2){\bf m}(a_2)k_2,{\bf t}(a){\bf n}(x_3){\bf m}(a_3)k_3)$ with $a, a_1,a_2,a_3 \in \R^\times$, $x_1,x_2,x_3\in\R$, and $k_1,k_2,k_3 \in {\rm SO}(2)$, where $da,da_i,dx_j$ are the Lebesgue measures and ${\rm vol}({\rm SO}(2),dk_i)=1$.

\subsubsection{Integral representation}\label{SS:integral rep.}

Let $\psi = \bigotimes_v\psi_v$ be a non-trivial additive character of $\F\backslash\A_\F$. 
For each place $v$ of $\F$, let $\mathcal{W}(\itPi_v,\psi_v\circ{\rm tr}_{\E_v/\F_v})$ be the space of Whittaker functions of $\itPi_v$ with respect to $\psi_v\circ{\rm tr}_{\E_v/\F_v}$.
When $v$ is finite, $\E_v/\F_v$ is unramified, $\itPi_v$ is unramified, and $\psi_v$ has conductor $\o_{\F_v}$, let $W_v^{\circ}$ be the $\GL_2(\o_{\E_v})$-invariant Whittaker function normalized so that $W_v^\circ(1)=1$.
For $\varphi \in V_\itPi$, let $W_{\varphi}$ be the Whittaker function of $\varphi$ with respect to $\psi\circ{\rm tr}_{\E/\F}$ defined by
\[
W_{\varphi,\psi}(g) = \int_{\E\backslash \A_\E}\varphi({\bf n}(x)g)\psi(-{\rm tr}_{\E/\F}(x))\,dx^{\rm Tam}.
\]
Here $dx^{\rm Tam}$ is the Tamagawa measure on $\E\backslash \A_\E$.
For automorphic forms $\varphi_1$ and $\varphi_2$ on $\G(\A_\F)$, assume one of them is cuspidal, we define the Petersson bilinear pairing
\[
\<\varphi_1,\varphi_2\>_{\G} = \int_{\A_\F^\times\G(\F)\backslash\G(\A_\F)}\varphi_1(g)\varphi_2(g)\,dg^{\rm Tam}.
\]
Here $dg^{\rm Tam}$ is the Tamagawa measure on $\A_\F^\times\backslash\G(\A_\F)$.

For $\varphi \in V_\itPi$ and holomorphic section $f^{(s)}$ of $I(\omega,s)$, by an unfolding argument (cf.\,\cite[\S\,2]{PSR1987}), we have
\[
\<\varphi,E(f^{(s)})\>_\G = |D_\F|^{-1/2}|D_\E|^{-1}\zeta_{\E}(2)^{-1}\cdot\int_{\A_\F^\times U^0(\A_\F)\backslash\G(\A_\F)}W_\varphi(g)f^{(s)}(\eta_\F\cdot g)\,dg
\]
for ${\rm Re}(s)$ sufficiently large.
Here $\zeta_\E(s)$ is the completed Dedekind zeta function of $\E$ and $\eta_\F \in G(\F)$ is any global element such that $Y(\F)\cdot\eta_\F = Y^0(\F)$.
We fix a choice of $\eta_\F$.
Note that the factor $|D_\F|^{-1/2}|D_\E|^{-1}\zeta_{\E}(2)^{-1}$ is due to the comparison between $dg^{\rm Tam}$ and $dg$ (cf.\,\cite[\S\,6.1]{IP2018}).
For each place $v$ of $\F$, $W_v \in \mathcal{W}(\itPi_v,\psi_v\circ{\rm tr}_{\E_v/\F_v})$, and $f_v^{(s)}$ a meromorphic section of $I(\omega_v,s)$, let $Z(W_v,f_v^{(s)})$ be the local zeta integral defined by
\begin{align}\label{E:local zeta}
Z(W_v,f_v^{(s)}) = \int_{\F_v^\times U^0(\F_v)\backslash \G(\F_v)}W_v(g_v)f_v^{(s)}(\eta_\F\cdot g_v)\,dg_v.
\end{align}
Note that the definition is not purely local, as it depends on $\eta_\F$. However, if we replace $\eta_\F$ by $\eta_v \in G(\F_v)$ such that $Y(\F_v)\cdot\eta_v = Y^0(\F_v)$, then the two local zeta integrals coincide up to a scalar determined by $\eta_\F\eta_v^{-1} \in P(\F_v)$.
Now we state the integral representation in the following proposition.
Let $S$ be a finite set of places of $\F$ including $S_\infty$ so that for $v \notin S$, 
\begin{itemize}
\item $\E_v/\F_v$ is unramified,
\item $\itPi_v$ is unramified, 
\item $\psi_v$ has conductor $\o_{\F_v}$,
\item $\nu(\eta_\F),\det (\eta_\F\vert_{Y(\F)}) \in \o_{\F_v}^\times$.
\end{itemize}
Here $\det (\eta\vert_{Y(\F)})$ is defined with respect to the isomorphism (\ref{E:iso}).
By \cite[Theorem 3.1]{PSR1987}, we have
\[
Z(W_v^\circ,f_{v,\circ}^{(s)}) = L(s+\tfrac{1}{2},\itPi_v,{\rm As})
\]
for all $v \notin S$. Thus we obtain the following

\begin{prop}[Garrett, Piatetski-Shapiro--Rallis]\label{P:integral rep.}
Let $W_S \in \mathcal{W}(\itPi_S,\psi_S\circ{\rm tr}_{\E_S/\F_S})$ and $f_S^{(s)}$ be a meromorphic section of $I(\omega_S,s)$.
Let $\varphi \in V_\itPi$ and $f^{(s)}$ the meromorphic section of $I(\omega,s)$ defined by
\[
W_{\varphi,\psi} = \prod_{v \notin S}W_v^\circ \cdot W_S,\quad f^{(s)} = \bigotimes_{v \notin S} f_{v,\circ}^{(s)}\otimes f_S^{(s)}.
\]
Then we have
\[
\<\varphi,E(f^{(s)})\>_\G = |D_\F|^{-1/2}|D_\E|^{-1}\zeta_{\E}(2)^{-1}\cdot L^S(s+\tfrac{1}{2},\itPi,{\rm As})\cdot Z(W_S,f_S^{(s)}).
\]
\end{prop}

\section{Proof of main results}\label{S:proofs}

\subsection{Proof of Theorem \ref{T:main 3}}\label{SS:proof 1}

We keep the notation of \S\,\ref{S:triple}. 
Let $\psi_\F = \bigotimes_v\psi_v$ be the standard additive character of $\F\backslash\A_\F$.
We assume further that $\itPi$ is cohomological. 
Let ${\sf w}$ be the integer such that 
\[
|\omega| = |\mbox{ }|_{\A_\F}^{\sf w}. 
\]
For $v \in S_\infty$, let $w_{1,v},w_{2,v},w_{3,v}$ be the places of $\E$ lying over $v$. We have
\[
\itPi_v = \boxtimes_{i=1}^3 D_{\kappa({w_{i,v}})} \otimes |\mbox{ }|_v^{{\sf w}(w_{i,v})/2}
\]
as representations of $({\rm R}_{\E_v/\F_v}\GL_2)(\F_v) = \GL_2(\R)^3$ for some $\kappa(w_{i,v}) \in \Z_{\geq 2}$ and ${\sf w}(w_{i,v}) \in \Z$ such that $\kappa(w_{i,v}) \equiv {\sf w}(w_{i,v}) \,({\rm mod}\,2)$. Here $D_\kappa$ is the discrete series representation of $\GL_2(\R)$ with weight $\kappa \in \Z_{\geq 2}$. We put $\underline{\kappa}_v = (\kappa(w_{1,v}),\kappa(w_{2,v}),\kappa(w_{3,v}))$ and $\underline{\kappa} = (\underline{\kappa}_v)_{v \in S_\infty}$.
For $\sigma \in {\rm Aut}(\C)$, define ${}^\sigma\!\underline{\kappa} = ({}^\sigma\!\underline{\kappa}_v)_{v \in S_\infty}$ with
\[
{}^\sigma\!\underline{\kappa}_v = (\kappa(\sigma^{-1}\circ w_{1,v}),\kappa(\sigma^{-1}\circ w_{2,v}),\kappa(\sigma^{-1}\circ w_{3,v})).
\]
Let $V_\itPi^-$ be the space of anti-holomorphic cusp forms in $V_\itPi$, that is, for $\varphi \in V_\itPi^-$ we have
\[
\varphi(g\cdot (k_{\theta_{1,v}},k_{\theta_{2,v}},k_{\theta_{3,v}})) = e^{-\sqrt{-1}\,(\theta_{1,v}\kappa(w_{1,v})+\theta_{2,v}\kappa(w_{2,v})+\theta_{3,v}\kappa(w_{3,v}))}\varphi(g)
\]
for all $g \in \GL_2(\A_\E)$ and all $v \in S_\infty$ and $(k_{\theta_{1,v}},k_{\theta_{2,v}},k_{\theta_{3,v}}) \in {\rm SO}(2)^3 \subset ({\rm R}_{\E_v/\F_v}\GL_2)(\F_v)$. Here $k_\theta = \bp \cos \theta & \sin\theta \\ -\sin\theta & \cos\theta\ep$.
For $v \in S_\infty$, let $W_{\underline{\kappa}_v}^- \in \mathcal{W}(\itPi_v,\psi_v\circ{\rm tr}_{\E_v/\F_v})$ be the Whittaker function of weight $-\underline{\kappa}_v$ normalized so that 
\[
W_{\underline{\kappa}_v}^-({\bf t}(-1)) = e^{-6\pi}.
\]
For $\varphi \in V_\itPi^-$, let $W_\varphi \in \prod_{v \nmid \infty}\mathcal{W}(\itPi_v,\psi_v\circ{\rm tr}_{\E_v/\F_v})$ be the unique Whittaker function on $\itPi_f = \bigotimes_{v \nmid \infty}\itPi_v$ such that
\begin{align}\label{E:anti-holo.}
W_{\varphi,\psi_\F} = \prod_{v \in S_\infty}W_{\underline{\kappa}_v}^-\cdot W_\varphi.
\end{align}
For $\sigma \in {\rm Aut}(\C)$, let ${}^\sigma\!\varphi \in V_{{}^\sigma\!\itPi}^-$ be the anti-holomorphic cusp form in $V_{{}^\sigma\!\itPi}$ defined so that
\[
W_{{}^\sigma\!\varphi,\psi_\F} = \prod_{v \in S_\infty}W_{{}^\sigma\!\underline{\kappa}_v}^-\cdot t_\sigma W_\varphi,
\]
where $t_\sigma W_\varphi(g_f) = \sigma(W({\bf t}(u_\sigma^{-1})g_f))$ for $g_f \in \GL_2(\A_{\E,f})$ and $u_\sigma \in \hat{\o}_{\F}^\times$ is the unique element such that $\sigma(\psi_\F(x)) = \psi_\F(u_\sigma x)$ for $x \in \A_{\F,f}$.

We assume the balanced condition (\ref{E:balanced 2}) is satisfied. Now we begin the proof of Theorem \ref{T:main 3}.
Let ${\rm Crit}^+(\itPi,{\rm As})$ be the set of right-half critical points of the triple product $L$-function $L(s,\itPi,{\rm As})$ given by
\begin{align*}
&{\rm Crit}^+(\itPi,{\rm As})\\
& = \left\{m+\tfrac{1}{2}\,\left\vert\,0\leq m+\tfrac{\sf w}{2} \leq \tfrac{\kappa(w_{1,v})+\kappa(w_{2,v})+\kappa(w_{3,v})-2\max\{\kappa(w_{1,v}),\kappa(w_{2,v}),\kappa(w_{3,v})\}}{2} -1 \mbox{ for all }v \in S_\infty\right\}\right..
\end{align*}
Firstly we consider the right-half critical points and explain how to reduce the problem to the construction of Eisenstein series in \S\,\ref{S:Eisenstein series} and computation of non-archimedean local zeta integrals in \S\,\ref{S:local zeta integral}.
For the left-half critical points, the proof will be given in \S\,\ref{SS:LHC} below.
Let $p(\itPi,{\rm As}) \in \C^\times$ be the period attached to the triple product $L$-function defined by
\[
p(\itPi,{\rm As}) = (|D_\E|/|D_\F|)^{1/2}\cdot(2\pi\sqrt{-1})^{\sum_{v \in S_\infty}(\kappa(w_{1,v})+\kappa(w_{2,v})+\kappa(w_{3,v})+2{\sf w}+2)}\cdot (\sqrt{-1})^{d{\sf w}}\cdot G(\omega)^2\cdot \Vert f_\itPi \Vert.
\]
For $1\leq i < j \leq 3$, let $X_{ij}$ be the weight raising differential operator in the complexified Lie algebra ${\rm Lie}(\GSp_6(\R))_\C = \frak{gsp}_6(\R)_\C$ defined by
\[
X_{ij} = \frac{1}{2\pi\sqrt{-1}}\bp\sqrt{-1}\,(e_{ij}+e_{ji}) & -(e_{ij}+e_{ji}) \\ -(e_{ij}+e_{ji}) & -\sqrt{-1}\,(e_{ij}+e_{ji})\ep.
\]
Here $e_{ij}$ is the $3$ by $3$ matrix with $(i,j)$-entry equals $1$ and zero otherwise. Then $X_{12}, X_{13},X_{23}$ take weight $(a,b,c)$ to weight $(a+1,b+1,c),(a+1,b,c+1),(a,b+1,c+1)$, respectively. For $\underline{\lambda} = (\lambda_1,\lambda_2,\lambda_3) \in \Z^3$ satisfying (\ref{E:balanced}) and $\ell\in\Z$ with 
\[
\ell \leq \lambda_1+\lambda_2+\lambda_3-2\max\{\lambda_1,\lambda_2,\lambda_3\}, \quad\ell \equiv \lambda_1+\lambda_2+\lambda_3 \,({\rm mod}\,2),
\]
let $X(\ell;\underline{\lambda}) \in \mathcal{U}(\frak{gsp}_6(\R)_\C)$ be the weight raising differential operator defined by
\[
X(\ell;\underline{\lambda}) = X_{12}^{\tfrac{\lambda_1+\lambda_2-\lambda_3-\ell}{2}}X_{13}^{\tfrac{\lambda_1+\lambda_3-\lambda_2-\ell}{2}}X_{23}^{\tfrac{\lambda_2+\lambda_3-\lambda_1-\ell}{2}}.
\]
Then it is clear that $X(\ell;\underline{\lambda})$ takes weight $(\ell,\ell,\ell)$ to weight $\underline{\lambda}$.
For $m+\tfrac{1}{2} \in {\rm Crit}^+(\itPi,{\rm As})$, let $\ell(m)$ be the integer defined by
\[
\ell(m) = 2m+{\sf w}+2.
\]
Then it is clear that the map $m+\tfrac{1}{2}\mapsto \ell(m)$ defines a bijective map from ${\rm Crit}^+(\itPi,{\rm As})$ to 
\begin{align*}
&S_{ad}\\
& = \left\{2 \leq \ell \leq \kappa(w_{1,v})+\kappa(w_{2,v})+\kappa(w_{3,v})-2\max\{\kappa(w_{1,v}),\kappa(w_{2,v}),\kappa(w_{3,v})\} \mbox{ for all $v \in S_\infty$},\,\ell \equiv {\sf w} \,({\rm mod}\,2)\right\}.
\end{align*}
Here $S_{ad}$ refers to the set of admissible weights for the Eisenstein series.
For $\ell \in S_{ad}$, let $X(\ell;\underline{\kappa})$ be the weight raising differential operator in the universal  enveloping algebra of the complexified Lie algebra of $({\rm R_{\E/\F}\GL_2})(\F_\infty)$ defined by
\[
X(\ell;\underline{\kappa}) = \bigotimes_{v \in S_\infty}X(\ell;\underline{\kappa}_v).
\]
We have the following result of Garrett--Harris in \cite[\S\,4]{GH1993}, which is proved base on the Serre duality for coherent cohomology of certain automorphic line bundles on the Shimura variety associated to $\G$ and the arithmeticity of $X(\ell;\underline{\kappa})$ (cf.\,\cite[\S\,2]{GH1993}).

\begin{thm}[Garrett--Harris]\label{T:GH1}
Let $\varphi \in V_\itPi^-$ and $\varphi'$ be a holomorphic automorphic form on $\GSp(V)(\A_\F)$ of weight $\ell \in S_{ad}$ and central character $\omega^{-1}$. Then we have
\[
\sigma \left( \frac{\<\varphi,\,X(\ell;\underline{\kappa})\cdot \varphi'\>_\G}{\Vert f_\itPi\Vert}\right) = \frac{\<{}^\sigma\!\varphi,\,X(\ell;{}^\sigma\!\underline{\kappa})\cdot {}^\sigma\!\varphi'\>_\G}{\Vert f_{{}^\sigma\!\itPi}\Vert}
\]
for all $\sigma \in {\rm Aut}(\C)$.
\end{thm}

Another key ingredient is the following result on the rationality of the archimedean local zeta integrals proved in \cite[\S\,5]{GH1993}. Note that the extra factor on the right-hand side is due to the comparison between $\eta_\F$ and $\theta$ in \cite[p.\,206]{GH1993}.

\begin{thm}[Garrett--Harris]\label{T:GH2}
Let $v \in S_\infty$ and $m+\tfrac{1}{2} \in {\rm Crit}^+(\itPi,{\rm As})$. We have
\begin{align*}
&\frac{Z(W_{\underline{\kappa}_v}^-,f_{v,\ell(m)}^{(s)}) \vert_{s=m}}{(2\pi\sqrt{-1})^{-2m+\kappa(w_{1,v})+\kappa(w_{2,v})+\kappa(w_{3,v})-{\sf w}-5}\cdot (\sqrt{-1})^{{\sf w}+1}}\in\iota_v(\nu(\eta_\F)^{-3m-2{\sf w}-3}\cdot\det(\eta_\F\vert_{Y(\F)})^{-2m-{\sf w}-2})\cdot \Q^\times.
\end{align*}
\end{thm}

Let $m+\tfrac{1}{2} \in {\rm Crit}^+(\itPi,{\rm As})$. Put $f_{\infty,\ell(m)}^{(s)} = \bigotimes_{v \in S_\infty}f_{v,\ell(m)}^{(s)} \in \bigotimes_{v\in S_\infty}I(\omega_v,s)$ (cf.\,(\ref{E:archimedean section}) and (\ref{E:degenerate p.s})). 
Let $S$ be a finite set of finite places satisfying conditions in \S\,\ref{SS:integral rep.} for $v \notin S\cup S_\infty$.
Let $W_S \in \mathcal{W}(\itPi_S,\psi_S\circ{\rm tr}_{\E_S/\F_S})$ and $f_S^{(s)}$ be a rational section of $I(\omega_S,s)$.
Let $\varphi \in V_\itPi^-$ be the anti-holomorphic cusp form defined so that (cf.\,(\ref{E:anti-holo.}))
\[
W_{\varphi} = \prod_{v \notin S\cup S_\infty}W_v^\circ \cdot W_S.
\]
By Proposition \ref{P:integral rep.} and Theorem \ref{T:GH2}, there exists $C \in \Q^\times$ such that
\begin{align}\label{E:Main thm 1 proof 1}
\begin{split}
&\left\<\varphi,X(\ell(m);\underline{\kappa})\cdot E\left(f_{\infty,\ell(m)}^{(s)}\otimes \bigotimes_{v \notin S\cup S_\infty} f_{v,\circ}^{(s)}\otimes f_S^{(s)}\right)\right\>_\G\\
& = C\cdot |D_\F|^{-1/2}|D_\E|^{-1}\zeta_{\E}(2)^{-1}\cdot (2\pi\sqrt{-1})^{2dm+\sum_{v \in S_\infty}(-\kappa(w_{1,v})-\kappa(w_{2,v})-\kappa(w_{3,v})+{\sf w}+5)}\cdot (\sqrt{-1})^{d{\sf w}+d}\\
&\times L^S(s+\tfrac{1}{2},\itPi,{\rm As})\cdot Z(W_S,f_S^{(s)}),\\
&\left\<{}^\sigma\!\varphi,X(\ell(m);{}^\sigma\!\underline{\kappa})\cdot E\left(f_{\infty,\ell(m)}^{(s)}\otimes \bigotimes_{v \notin S\cup S_\infty} {}^\sigma\!f_{v,\circ}^{(s)}\otimes {}^\sigma\!f_S^{(s)}\right)\right\>_\G\\
& = C\cdot |D_\F|^{-1/2}|D_\E|^{-1}\zeta_{\E}(2)^{-1}\cdot (2\pi\sqrt{-1})^{2dm+\sum_{v \in S_\infty}(-\kappa(w_{1,v})-\kappa(w_{2,v})-\kappa(w_{3,v})+{\sf w}+5)}\cdot (\sqrt{-1})^{d{\sf w}+d}\\
&\times L^S(s+\tfrac{1}{2},{}^\sigma\!\itPi,{\rm As})\cdot Z(t_{\sigma,S}W_S,{}^\sigma\!f_S^{(s)}).
\end{split}
\end{align}
Note that we have incorporate ${\rm N}_{\F/\Q}(\nu(\eta_\F)^{-3m-2{\sf w}-3}\cdot\det(\eta_\F\vert_{Y(\F)})^{-2m-{\sf w}-2})$ into the constant $C$. 
Now we assume the rational section $f_S^{(s)}$ is chosen so that the following conditions are satisfied:
\begin{itemize}
\item[(1)] The Eisenstein series $E\left(f_{\infty,\ell(m)}^{(s)}\otimes \bigotimes_{v \notin S \cup S_\infty} f_{v,\circ}^{(s)}\otimes f_S^{(s)}\right)$ is holomorphic at $s=m$.
\item[(2)] $E^{[\ell(m)]}\left(\bigotimes_{v \notin S \cup S_\infty} f_{v,\circ}^{(s)}\otimes f_S^{(s)}\right)$ is a holomorphic automorphic form of weight $\ell(m)$ and satisfies (\ref{E:Galois equiv. E.S.}) for $n=3$.
\end{itemize}
Under the assumption, we then deduce from Theorem \ref{T:GH1} and (\ref{E:Main thm 1 proof 1}) that
\begin{align}\label{E:Main thm 1 proof 2}
\sigma \left( \frac{\left.\left(L^S(s+\tfrac{1}{2},\itPi,{\rm As})\cdot Z(W_S,f_S^{(S)})\right)\right\vert_{s=m}}{(2\pi\sqrt{-1})^{4dm}\cdot p(\itPi,{\rm As})}\right) = \frac{\left.\left(L^S(s+\tfrac{1}{2},{}^\sigma\!\itPi,{\rm As})\cdot Z(t_{\sigma,S}W_S,{}^\sigma\!f_S^{(S)})\right)\right\vert_{s=m}}{(2\pi\sqrt{-1})^{4dm}\cdot p({}^\sigma\!\itPi,{\rm As})}
\end{align}
for all $\sigma \in{\rm Aut}(\C)$.
Here we have used the fact that $\zeta_\E(2) \in |D_\E|^{1/2}\cdot\pi^3\cdot\Q^\times$.
We also have the following result on the Galois equivariant property of the non-archimedean local zeta integrals.

\begin{lemma}\label{L:Galois equiv. local zeta}
Let $v$ be a finite place. Let $W_v \in \mathcal{W}(\itPi_v,\psi_v\circ{\rm tr}_{\E_v/\F_v})$ and $f_v^{(s)}$ be a rational section of $I(\omega_v,s)$. 
Then we have
\[
{}^\sigma\!Z(W_v,f_v^{(s)}) = Z(t_{\sigma,v}W_v,{}^\sigma\!f_v^{(s)})
\]
as rational functions in $q_v^{-s}$ for all $\sigma \in {\rm Aut}(\C)$.
\end{lemma}

\begin{proof}
The assertion can be proved following arguments similar to \cite[Proposition 5.14]{Chen2021c}.
\end{proof}

We assume further that the Whittaker function $W_S$ is chosen so that the following condition is satisfied:
\begin{itemize}
\item[(3)] The local zeta integral $Z(W_S,f_S^{(s)})$ is holomorphic and non-vanishing at $s = m$.
\end{itemize}
Then it follows from (\ref{E:Main thm 1 proof 2}) and Lemma \ref{L:Galois equiv. local zeta} that 
\[
\sigma \left( \frac{L^S(m+\tfrac{1}{2},\itPi,{\rm As})}{(2\pi\sqrt{-1})^{4dm}\cdot p(\itPi,{\rm As})}\right) = \frac{L^S(m+\tfrac{1}{2},{}^\sigma\!\itPi,{\rm As})}{(2\pi\sqrt{-1})^{4dm}\cdot p({}^\sigma\!\itPi,{\rm As})}
\]
for all $\sigma \in{\rm Aut}(\C)$.
The algebraicity for $L^{(\infty)}(m+\tfrac{1}{2},\itPi,{\rm As})$ then follows from the Galois equivariant property (\ref{E:local factor non-archimedean}) of the local factors and the following lemma.

\begin{lemma}
Let $v$ be a finite place. Then the local factor $L(s+\tfrac{1}{2},\itPi,{\rm As})$ is holomorphic for ${\rm Re}(s)\geq -\tfrac{\sf w}{2}$.
\end{lemma}

\begin{proof}
The assertion is a direct consequence of the result of Kim and Shahidi \cite{KS2002} and \cite[Lemma 3.1]{Chen2021b}.
\end{proof}

In conclusion, in order to prove Theorem \ref{T:main 3} for $m+\tfrac{1}{2} \in {\rm Crit}^+(\itPi,{\rm As})$, we are reduced to show that conditions (1)-(3) are satisfied for some $f_S^{(s)}$ and $W_S$ when $2 \leq \ell(m) \leq 4$.
Indeed, by Theorem \ref{T:Harris}, conditions (1) and (2) are always satisfied when $\ell(m) >4$ for any rational section $f_S^{(s)}$ holomorphic for ${\rm Re}(s) \geq 3-\tfrac{\sf w}{2}$.
Once we have freedom to choose $f_S^{(s)}$, condition (3) is also satisfied by \cite[Proposition 3.3]{PSR1987}.
For conditions (1) and (2), we have proved in Proposition \ref{P:holomorphy} how to construct holomorphic Eisenstein series satisfying (\ref{E:Galois equiv. E.S.}) when $2 \leq \ell(m) \leq 4$.
Invoking the assumptions in Proposition \ref{P:holomorphy} for $n=3$, to verify conditions (1)-(3), we are further reduced to verify the following condition:
\begin{align}\label{E:condition}
\begin{split}
&\mbox{There exist a finite place $v$, $\itPhi_v \in \mathcal{S}({\rm Sym}_3(\F_v))$, and a Whittaker function $W_v$ of $\itPi_v$ such that}\\
&\mbox{$\widehat{\Phi}_v$ is supported in ${\rm Sym}_3(\F_v) \cap \GL_3(\F_v)$ and $Z(W_v,f_{\itPhi_v}^{(s)})$ is holomorphic and non-vanishing at $s=m$}.
\end{split}
\end{align}
For a finite place $v$ splits in $\E$ and such that $\itPi_v$ is unramified, we write $\itPi_v = \itPi_{v}^{(1)}\boxtimes \itPi_{v}^{(2)} \boxtimes \itPi_{v}^{(3)}$ for some irreducible generic unramified representation $\itPi_{v}^{(i)}$ of $\GL_2(\F_v)$. For $i=1,2,3$, let $\alpha_{i,v},\beta_{i,v}$ be the Satake parameters of $\itPi_{i,v}$.
By the result of Kim and Shahidi \cite{KS2002}, we have
\begin{align}\label{E:KS estimation}
q_v^{-1/2} < q_v^{{\sf w}/2}|a_{1,v}a_{2,v}a_{3,v}| < q_v^{1/2}
\end{align}
for $a_{i,v} \in \{\alpha_{i,v},\beta_{i,v}\}$.
We have the following result on the explicit computation of non-archimedean local zeta integrals. See also Remark \ref{R:support}.
\begin{prop}\label{P:local zeta} 
Let $v\nmid 2$ be a finite place $v$ splits in $\E$ and such that $\itPi_v$ is unramified. Then there exist $\itPhi_v,\itPhi_v' \in \mathcal{S}({\rm Sym}_3(\F_v))$ and Whittaker functions $W_v,W_v'$ of $\itPi_v$ such that $\widehat{\Phi}_v$ and $\widehat{\Phi}_v'$ are supported in ${\rm Sym}_3(\F_v) \cap \GL_3(\F_v)$ and
\begin{align*}
&Z(W_v,f_{\itPhi_v}^{(s)})\vert_{s=m} \\
& = 
(q_v^{m-1/2}-\alpha_{1,v}\beta_{2,v}\beta_{3,v})(q_v^{m-1/2}-\beta_{1,v}\beta_{2,v}\beta_{3,v})(q_v^{m-1/2}-\alpha_{1,v}\alpha_{2,v}\beta_{3,v})(q_v^{m-1/2}-\beta_{1,v}\alpha_{2,v}\beta_{3,v}),\\
&Z(W_v',f_{\itPhi_v'}^{(s)})\vert_{s=m}\\
& = (q_v^{m-1/2}-\alpha_{1,v}\beta_{2,v}\beta_{3,v})(q_v^{m-1/2}-\beta_{1,v}\beta_{2,v}\beta_{3,v})(q_v^{m-1/2}-\beta_{1,v}\alpha_{2,v}\beta_{3,v})(q_v^{m-1/2}-\beta_{1,v}\beta_{2,v}\alpha_{3,v}).
\end{align*}
\end{prop}

\begin{proof}
The assertion will be proved in Lemma \ref{L:local zeta 3} below. Note that we have incorporate some fudge factors into the definition of $\itPhi_v$ and $\itPhi_v'$.
\end{proof}

When $\ell(m)\in\{2,4\}$, we have $m \in\{ -\tfrac{\sf w}{2},-\tfrac{\sf w}{2}+1\}$.
In these cases, we see that condition (\ref{E:condition}) holds by estimation (\ref{E:KS estimation}) and Proposition \ref{P:local zeta}.
In the rest of the proof, we deal with the remaining case $\ell(m)=3$, that is, $m=-\tfrac{\sf w}{2}+\tfrac{1}{2}$.
In this case, ${\sf w}$ must be odd and hence the character $\omega|\mbox{ }|_{\A_{\F}}^{-{\sf w}}$ is non-trivial as $\omega_v(-1) = (-1)^{\sf w}$ for all $v \in S_\infty$.
Let $S_{good}$ be the set of finite places of $\F$ defined by
\begin{align}\label{E:good places}
S_{good} = \left\{ v \nmid2\cdot\infty\,\vert\, \mbox{$v$ splits in $\E$, $\itPi_v$ is unramified, and $\omega_v\neq|\mbox{ }|_v^{{\sf w}}$}\right\}.
\end{align}

\begin{lemma}\label{L:density}
The set $S_{good}$ has positive density.
\end{lemma}

\begin{proof}
Let $\L_\omega$ be the ray class field associated to $\omega$.
The assertion is clear if $\E = \F \times \F \times \F$.
If $\E = \F'\times\F$ for some totally real quadratic extension $\F'/\F$, consider $\F'\cap\L_\omega$.
When $\F'\cap\L_\omega = \F$, the assertion follows from Chebotarev's density theorem.
When $\F' \subset \L_\omega$, $S_{good}$ must has positive density. Indeed, if otherwise, then $\omega = |\mbox{ }|_{\A_\F}^{\sf w}\omega_{\F'/\F}$ by Chebotarev's density theorem for $\L_\omega/\F$. Here $\omega_{\F'/\F}$ is the quadratic Hecke character of $\A_\F^\times$ associated to $\F'/\F$ by class field theory. This contradicts the oddness of ${\sf w}$.
If $\E$ is a field, consider $\widetilde{\E}\cap\L_\omega$, where $\widetilde{\E}$ is the Galois closure of $\E/\F$. When $\widetilde{\E}\cap\L_\omega = \F$, the assertion follows from Chebotarev's density theorem for $\widetilde{\E}\L_\omega/\F$. When $\widetilde{\E}\cap\L_\omega  \neq \F$ and $\E/\F$ is Galois, we have $\E \subset \L_\omega$ and similarly as above $S_{good}$ has positive density. If $\widetilde{\E}\cap\L_\omega \neq \F$ and $\E/\F$ is not Galois, then $\widetilde{\E}\cap\L_\omega$ is the unique totally real quadratic extension over $\F$ inside $\widetilde{\E}$. Indeed, $\widetilde{\E}\cap\L_\omega / \F$ is abelian whereas $\widetilde{\E}/\F$ is totally real and non-abelian. In this case, similarly $S_{good}$ has positive density by applying Chebotarev's density theorem to $\widetilde{\E}\L_\omega / \F$ together with the oddness of ${\sf w}$.
This completes the proof.
\end{proof}

\begin{lemma}\label{L:good places}
Let $v \in S_{good}$. Then either one of the following cases hold:
\begin{itemize}
\item[(1)] Up to permutations on $\{\itPi_{v}^{(1)},\itPi_{v}^{(2)},\itPi_{v}^{(3)}\}$ and on the Satake parameters of $\itPi_{v}^{(i)}$, we have
\[
(q_v^{-{\sf w}/2}-\alpha_{1,v}\beta_{2,v}\beta_{3,v})(q_v^{-{\sf w}/2}-\beta_{1,v}\beta_{2,v}\beta_{3,v})\neq0
\]
\item[(2)] We have
\[
\alpha_{i,v}+\beta_{i,v}=0,\quad\omega_v^2=|\mbox{ }|_v^{2{\sf w}}
\]
for $i=1,2,3$.
\end{itemize}
\end{lemma}

\begin{proof}
We drop the subscript $v$ for brevity.
If $a_1a_2a_3 \neq q^{-{\sf w}/2}$ for all $a_i \in \{\alpha_i,\beta_i\}$ for $i=1,2,3$, then clearly (1) holds.
Suppose otherwise and $\alpha_1\alpha_{2}\alpha_{3} = q^{-{\sf w}/2}$. Then $\omega_v\neq|\mbox{ }|_v^{{\sf w}}$ implies that $\beta_{1}\beta_{2}\beta_{3} \neq q^{-{\sf w}/2}$.
In this case, consider $\alpha_{1}\beta_{2}\beta_{3}$, $\beta_{1}\alpha_{2}\beta_{3}$, and $\beta_{1}\beta_{2}\alpha_{3}$.
If any one of them is not equal to $q^{-{\sf w}/2}$, then we are in Case (1). Otherwise we would have 
\begin{align}\label{E:good places proof 1}
\alpha_{1}\beta_{2}\beta_{3}= q^{-{\sf w}/2},\quad \beta_{1}\alpha_{2}\beta_{3}= q^{-{\sf w}/2},\quad \beta_{1}\beta_{2}\alpha_{3}= q^{-{\sf w}/2}.
\end{align}
Thus $\alpha_{1}\alpha_{2}\alpha_{3}(\beta_{1}\beta_{2}\beta_{3})^2 = q^{-3{\sf w}/2}$, which implies that $\beta_{1}\beta_{2}\beta_{3} = -q^{-{\sf w}/2}$.
In particular, we have $\omega(\varpi) = -q^{-{\sf w}}$.
Now we repeat the above argument with $\alpha_{1}\alpha_{2}\alpha_{3}$ replaced by $\alpha_{1}\beta_{2}\beta_{3}$, $\beta_{1}\alpha_{2}\beta_{3}$, and $\beta_{1}\beta_{2}\alpha_{3}$. Then we ended up with either (1) holds or 
\begin{align}\label{E:good places proof 2}
\beta_{1}\alpha_{2}\alpha_{3} = -q^{-{\sf w}/2},\quad \alpha_{1}\beta_{2}\alpha_{3} = -q^{-{\sf w}/2},\quad
\alpha_{1}\alpha_{2}\beta_{3}  = -q^{-{\sf w}/2}.
\end{align}
(\ref{E:good places proof 1}) and (\ref{E:good places proof 2}) then imply that $\alpha_{i}+\beta_{i}=0$ for $i=1,2,3$.
This completes the proof.
\end{proof}

\begin{lemma}\label{L:good places 2}
Let $v \in S_{good}$ and assume (1) of Lemma \ref{L:good places} holds. 
Then one of the local zeta integrals in Proposition \ref{P:local zeta} is non-zero.
\end{lemma}

\begin{proof}
We drop the subscript $v$ for brevity. Assume 
\[
(q^{-{\sf w}/2}-\alpha_{1}\beta_{2}\beta_{3})(q^{-{\sf w}/2}-\beta_{1}\beta_{2}\beta_{3})\neq0
\]
To prove the assertion, it suffices to show that either
\begin{align}\label{E:good places 2 proof 1}
(q^{-{\sf w}/2}-\alpha_{1}\gamma)(q^{-{\sf w}/2}-\beta_{1}\gamma) \neq 0
\end{align}
for some $\gamma \in \{\alpha_2\beta_3,\beta_2\alpha_3\}$ or
\begin{align}\label{E:good places 2 proof 2}
(q^{-{\sf w}/2}-\gamma\alpha_2\beta_3)(q^{-{\sf w}/2}-\gamma\beta_2\alpha_3) \neq 0
\end{align}
for some $\gamma \in \{\alpha_1,\beta_1\}$.
If $\alpha_1\alpha_2\beta_3=\beta_1\alpha_2\beta_3=q^{-{\sf w}/2}$, then the condition $\omega \neq |\mbox{ }|^{\sf w}$ implies that (\ref{E:good places 2 proof 1}) holds for $\gamma =\beta_2\alpha_3$.
If $\alpha_1\alpha_2\beta_3=q^{-{\sf w}/2}$ and $\beta_1\alpha_2\beta_3 \neq q^{-{\sf w}/2}$, then then the condition $\omega \neq |\mbox{ }|^{\sf w}$ implies that $\beta_1\beta_2\alpha_3 \neq q^{-{\sf w}/2}$. Thus (\ref{E:good places 2 proof 2}) holds for $\gamma=\beta_1$.
If $\alpha_1\alpha_2\beta_3\neq q^{-{\sf w}/2}$ and $\beta_1\alpha_2\beta_3 = q^{-{\sf w}/2}$, then similarly (\ref{E:good places 2 proof 2}) holds for $\gamma=\alpha_1$.
If $\alpha_1\alpha_2\beta_3\neq q^{-{\sf w}/2}$ and $\beta_1\alpha_2\beta_3 \neq q^{-{\sf w}/2}$, this is (\ref{E:good places 2 proof 1}) for $\gamma=\alpha_2\beta_3$.
This completes the proof.
\end{proof}

By Lemma \ref{L:good places 2}, we see that condition (\ref{E:condition}) holds if (1) of Lemma \ref{L:good places} holds for some $v \in S_{good}$.
In \S\,\ref{SS:Case 1}-\ref{SS:Case 3} below, we will prove the existence of such $v \in S_{good}$ except for one exceptional case. In the exceptional case, we need the following result on the algebraicity of ratios of critical $L$-values for algebraic Hecke characters of CM-fields.

\begin{lemma}\label{L:period relation}
Let $\chi$ be an algebraic Hecke character of $\A_{\mathbb{K}}^\times$ for some CM-extension $\mathbb{K}/\F$.
Assume $\chi_v$ is not equal to a power of $|\mbox{ }|_v\circ{\rm N}_{\mathbb{K}_v/\F_v}$ for all $v \in S_\infty$. Then for any finite order Hecke character $\eta$ of $\A_\F^\times$ with parallel signature $1$ and any critical point $m \in \Z$ of the Hecke $L$-function $L(s,\chi)$, we have
\[
\sigma\left( \frac{L^{(\infty)}(m,\chi\cdot(\eta\circ{\rm N}_{\KK/\F}))}{G(\mu)\cdot L^{(\infty)}(m,\chi)}\right) = \frac{L^{(\infty)}(m,{}^\sigma\!\chi\cdot({}^\sigma\!\eta\circ{\rm N}_{\KK/\F}))}{G({}^\sigma\!\mu)\cdot L^{(\infty)}(m,{}^\sigma\!\chi)}
\]
for all $\sigma \in {\rm Aut}(\C)$.
\end{lemma}

\begin{proof}
Let $I_{\mathbb{K}}^\F(\chi)$ be the automorphic induction of $\chi$ to $\GL_2(\A_\F)$. By the assumption on $\chi$, we know that $I_{\mathbb{K}}^\F(\chi)\otimes |\mbox{ }|_{\A_\F}^{-1/2}$ is a cohomological irreducible cuspidal automorphic representation of $\GL_2(\A_\F)$.
The assertion thus follows from the result of Shimura \cite[Theorem 4.3]{Shimura1978} (cf.\,Theorem \ref{T:Shimura} below) on the algebraicity of critical values of standard $L$-function for $\GL_2(\A_\F)$ and the equality of $L$-functions
\[
L(s,I_{\mathbb{K}}^\F(\chi)) = L(s,\chi).
\]
This completes the proof.
\end{proof}

We have the following three cases
\begin{align*}
\begin{cases}
\E=\F\times \F \times \F & \underline{\mbox{Case 1}},\\
\E=\F'\times \F \mbox{ for some real quadratic extension $\F'$ of $\F$}& \underline{\mbox{Case 2}},\\
\E \mbox{ is a field} & \underline{\mbox{Case 3}}.
\end{cases}
\end{align*}
Denote by $\L_\omega$ the ray class field associated to $\omega$ by class field theory.
Note that $\L_\omega$ is totally imaginary over $\F$ by our assumption that ${\sf w}$ is odd.
For the exceptional case, we also need the period relation
\begin{align}\label{E:Case 1 proof 1}
\sigma \left( \frac{p(\itPi\otimes\eta,{\rm As})}{G(\eta\vert_{\A_\F^\times})^4\cdot p(\itPi,{\rm As})}\right) = \frac{p({}^\sigma\!\itPi\otimes{}^\sigma\!\eta,{\rm As})}{G({}^\sigma\!\eta\vert_{\A_\F^\times})^4\cdot p({}^\sigma\!\itPi,{\rm As})}
\end{align}
for all $\sigma \in {\rm Aut}(\C)$. To verify it, note that $\omega_{\itPi\otimes\eta} = \omega_\itPi\eta^2$. In general $f_\itPi$ and $f_{\itPi \otimes \eta}$ are not equal. Nonetheless, we have
\[
\sigma \left( \frac{\Vert f_{\itPi\otimes\eta} \Vert}{\Vert f_\itPi \Vert}\right) = \frac{\Vert f_{{}^\sigma\!\itPi\otimes{}^\sigma\!\eta} \Vert}{\Vert f_{{}^\sigma\!\itPi} \Vert}
\]
for all $\sigma \in {\rm Aut}(\C)$. Thus (\ref{E:Case 1 proof 1}) holds.
In \S\,\ref{SS:Case 1}-\ref{SS:Case 3} below, we will repeatedly use Chebotarev's density theorem (CDT for abbreviation) and the following corollary to the Sato--Tate conjecture.
The conjecture in the CM case is a consequence of the equidistribution of eigenvalues of Hecke characters (cf.\,\cite[Chapter I, Appendix A.2, Theorem 2]{Serre1998}). For the non-CM case, it was proved by Barnet-Lamb, Gee, and Geraghty \cite{BLGG2011} (see also \cite{BLGHT2011} for the elliptic case).

\begin{lemma}\label{L:ST}
Let $\L$ be a totally real number field and $\itSigma$ a cohomological irreducible cuspidal automorphic representation of $\GL_2(\A_\L)$. For a finite place $v$ of $\L$ such that $\itSigma_v$ is unramified, denote by $C(\itSigma_v)\subset\GL_2(\C)$ the semi-simple conjugacy class associated to $\itSigma_v$. Let $S$ be a set of finite places of $\L$. 
\begin{itemize}
\item[(1)] If $\itSigma$ is non-CM and $S$ has positive density, then there exists $v \in S$ such that ${\rm tr}(C(\itSigma_v)) \neq 0$.
\item[(2)] If $\itSigma$ is CM by a CM-extension $\K/\L$ and $S\cap S_{\K}$ has positive density, then there exists $v \in S\cap S_{\K}$ such that ${\rm tr}(C(\itSigma_v)) \neq 0$.
Here $S_{\K}$ is the set of places of $\L$ that splits in $\K$.
\end{itemize}
Here by density we refer to the natural density.
\end{lemma}

\subsubsection{\underline{Case $1$}}\label{SS:Case 1}

In this case, we have $\itPi = \itPi_1 \boxtimes \itPi_2 \boxtimes \itPi_3$ for some cohomological irreducible cuspidal automorphic representation $\itPi_i$ of $\GL_2(\A_\F)$. 
We consider the two subcases:
\begin{itemize}
\item[\underline{1-(i)}] One of $\itPi_1,\itPi_2,\itPi_3$ is non-CM.
\item[\underline{1-(ii)}] For $i=1,2,3$, there exist CM-extension $\KK_i/\F$ and algebraic Hecke character $\chi_i$ of $\A_\KK^\times$ such that $\itPi_i = I_{\KK_i}^\F(\chi_i)\otimes|\mbox{ }|_{\A_\F}^{-1/2}$.
\end{itemize}
In Case 1-(i), assume $\itPi_1$ is non-CM. Since $S_{good}$ has non-zero density by Lemma \ref{L:density}, it follows from Lemma \ref{L:ST} for $\itPi_1$ that there exists $v \in S_{good}$ such that ${\rm tr}(C(\itPi_{1,v})) \neq 0$. 
Therefore (1) of Lemma \ref{L:good places} must holds for this $v$.

In Case 1-(ii), let $S_{\KK_i}$ be the set of places of $\F$ that splits in $\KK_i$ for $i=1,2,3$. 
\begin{itemize}
\item If $S_{good}\cap S_{\KK_i}$ has positive density for some $1 \leq i \leq 3$, then there exists $v \in S_{good}\cap S_{\KK_i}$ such that ${\rm tr}(C(\itPi_{i,v}))\neq0$ by Lemma \ref{L:ST} for $\itPi_i$. 
Therefore (1) of Lemma \ref{L:good places} holds for this $v$.
\item
If there exists $1 \leq i \leq 3$ such that $\K_i\cap\L_\omega = \F$, then $S_{good} \cap S_{\K_i}$ has positive density by CDT for $\KK_i\L_\omega$. Thus similarly (1) of Lemma \ref{L:good places} holds for some $v\in S_{good} \cap S_{\K_i}$.
\item
If $\KK_i\subset\L_\omega$ and $S_{good}\cap S_{\KK_i}$ has density zero for $i=1,2,3$, then $\omega = \omega_{\KK_i/\F}|\mbox{ }|_{\A_\F}^{\sf w}$ by applying CDT for $\L_\omega/\F$. Here $\omega_{\KK_i/\F}$ is the quadratic Hecke character of $\A_\F^\times$ associated to $\KK_i/\F$ by class field theory for $i=1,2,3$.
In particular, $\KK_1=\KK_2=\KK_3$. We write $\K$ for this common field. 
Then we have the factorization of $L$-functions:
\[
L^S(s+\tfrac{1}{2},\itPi,{\rm As}) = L^S(s-1,\chi_1\chi_2\chi_3)L^S(s-1,\chi_1^c\chi_2\chi_3)L^S(s-1,\chi_1\chi_2^c\chi_3)L^S(s-1,\chi_1\chi_2\chi_3^c)
\]
for sufficiently large finite set $S$ of places.
Here $c$ is the non-trivial automorphism in $\Gal(\KK/\F)$.
Let $\eta$ be a finite order Hecke character of $\A_\F^\times$. Then similarly we have
\begin{align*}
L^S(s+\tfrac{1}{2},\itPi\otimes\eta,{\rm As}) &= L^S(s-1,\chi_1\chi_2\chi_3\cdot(\eta\vert_{\A_\F^\times}\circ{\rm N}_{\K/\F}))L^S(s-1,\chi_1^c\chi_2\chi_3\cdot(\eta\vert_{\A_\F^\times}\circ{\rm N}_{\K/\F}))\\
&\times L^S(s-1,\chi_1\chi_2^c\chi_3\cdot(\eta\vert_{\A_\F^\times}\circ{\rm N}_{\K/\F}))L^S(s-1,\chi_1\chi_2\chi_3^c\cdot(\eta\vert_{\A_\F^\times}\circ{\rm N}_{\K/\F})).
\end{align*}
Assume $\eta\vert_{\A_\F^\times}$ has parallel signature $1$.
The unbalanced condition (\ref{E:balanced 2}) implies that the characters $\chi_1\chi_2\chi_3$, $\chi_1^c\chi_2\chi_3$, $\chi_1\chi_2^c\chi_3$, $\chi_1\chi_2\chi_3^c$ all satisfy the assumption in Lemma \ref{L:period relation}.
Thus $L(1-\tfrac{\sf w}{2},\itPi,{\rm As}) \neq 0$ and we have
\begin{align}\label{E:Case 1 proof 2}
\sigma \left( \frac{L^{(\infty)}(1-\tfrac{\sf w}{2},\itPi\otimes\eta,{\rm As})}{G(\eta\vert_{\A_\F^\times})^4\cdot L^{(\infty)}(1-\tfrac{\sf w}{2},\itPi,{\rm As})}\right) = \frac{L^{(\infty)}(1-\tfrac{\sf w}{2},{}^\sigma\!\itPi\otimes{}^\sigma\!\eta,{\rm As})}{G({}^\sigma\!\eta\vert_{\A_\F^\times})^4\cdot L^{(\infty)}(1-\tfrac{\sf w}{2},{}^\sigma\!\itPi,{\rm As})}
\end{align}
for all $\sigma \in {\rm Aut}(\C)$.
By (\ref{E:Case 1 proof 1}) and (\ref{E:Case 1 proof 2}), the algebraicity holds for $L^{(\infty)}(1-\tfrac{\sf w}{2},\itPi,{\rm As})$ if and only if it holds for $L^{(\infty)}(1-\tfrac{\sf w}{2},\itPi\otimes\eta,{\rm As})$.
Assume further $\eta$ is chosen so that $\eta\vert_{\A_\F^\times}$ has order $3$. Let $\L/\F$ be the cubic Galois extension associated to $\eta\vert_{\A_\F^\times}$ by class field theory. Then $\K\cap\L=\F$ and the following set of places of $\F$ has positive density by CDT for $\K\L/\F$:
\[
\left\{v \nmid 2\cdot\infty \,\left\vert\, \mbox{$v$ splits in $\KK$, $\itPi_v$ and $\eta_v$ are unramified, and $\eta\vert_{\F_v^\times} \neq 1$}\right\}\right..
\]
It then follows from Lemma \ref{L:ST} for $\itPi_1$ that there exists $v$ belongs to the above set such that ${\rm tr}(C(\itPi_{1,v})) \neq 0$. On the other hands, it is clear that the above set is contained in the set (\ref{E:good places}) with $\itPi$ and $\omega$ replaced by $\itPi\otimes \eta$ and $\omega\eta\vert_{\A_\F^\times}^2$. 
Indeed, since $v$ splits in $\KK$ and $\eta\vert_{\F_v^\times}$ has order $3$, we have $\omega_v\eta\vert_{\F_v^\times}^2 = |\mbox{ }|^{\sf w}_v\eta\vert_{\F_v^\times}^2 \neq |\mbox{ }|_v^{\sf w}$.
Thus (1) of Lemma \ref{L:good places} holds with $\itPi$ replaced by $\itPi\otimes\eta$, and the algebraicity holds for $L^{(\infty)}(1-\tfrac{\sf w}{2},\itPi\otimes\eta,{\rm As})$.
\end{itemize}

\subsubsection{\underline{Case $2$}}\label{SS:Case 2}

In this case, we have $\itPi = \itPi_1 \boxtimes \itPi_2$ for some cohomological irreducible cuspidal automorphic representations $\itPi_1$ and $\itPi_2$ of $\GL_2(\A_{\F'})$ and $\GL_2(\A_\F)$. 
We consider the two subcases:
\begin{itemize}
\item[\underline{2-(i)}] One of $\itPi_1$ and $\itPi_2$ is non-CM.
\item[\underline{2-(ii)}] There exist CM-extensions $\KK_1/\F'$ and $\KK_2/\F$, and algebraic Hecke characters $\chi_1$ of $\A_{\KK_1}^\times$ and $\chi_2$ of $\A_{\KK_2}^\times$ such that $\itPi_1 = I_{\KK_1}^{\F'}(\chi_1)\otimes|\mbox{ }|_{\A_{\F'}}^{-1/2}$ and $\itPi_2 = I_{\KK_2}^\F(\chi_2)\otimes|\mbox{ }|_{\A_{\F}}^{-1/2}$.
\end{itemize}
The proof for Case $2$-(i) is similar to that of Case 1-(i) by applying Lemmas \ref{L:density} and \ref{L:ST}.

In Case $2$-(ii), let $S_{\KK_i}$ be the set of places of $\F$ that splits in $\KK_i$ for $i=1,2$.
\begin{itemize}
\item
If $S_{good}\cap S_{\K_j}$ has positive density for some $j=1,2$, then there exists $v\in S_{good}\cap S_{\K_j}$ such that ${\rm tr}(C(\itPi_{v}^{(i)})) \neq 0$ for some $i=1,2,3$ by Lemma \ref{L:ST} for $\itPi_j$.
Therefore (1) of Lemma \ref{L:good places} holds for this $v$.
\item
If $\K_1/\F$ is cyclic, let $\omega_{\K_1/\F}$ and $\omega_{\F'/\F}$ be Hecke characters of order $4$ and $2$ associated to $\K_1/\F$ and $\F'/\F$, respectively, by class field theory.
Consider $\K_1\cap\L_\omega$. When $\K_1\cap\L_\omega= \F$, the set $S_{good}\cap S_{\K_1}$ has positive density by CDT for $\K_1\L_\omega / \F$.
Thus (1) of Lemma \ref{L:good places} holds for some $v \in S_{good}\cap S_{\K_1}$.
If $\K_1\cap\L_\omega  \neq \F$ and $S_{good}\cap S_{\K_1}$ has density zero, then $\omega$ is equal to either $\omega_{\K_1/\F}|\mbox{ }|_{\A_\F}^{\sf w}$ or $\omega_{\F'/\F}|\mbox{ }|_{\A_\F}^{\sf w}$.
This follows from applying CDT to $\K_1\L_\omega/\F$.
The case $\omega = \omega_{\F'/\F}|\mbox{ }|_{\A_\F}^{\sf w}$ is not possible by the oddness of ${\sf w}$. 
Since $\K_1/\F$ is cyclic and $\F'/\F$ is totally real, we have $\K_1\cap \K_2  = \F$. By CDT for $\K_1\K_2/\F$ and the condition that $\omega = \omega_{\K_1/\F}|\mbox{ }|_{\A_\F}^{\sf w}$, the set $S_{good}\cap S_{\K_2}$ has positive density. Thus (1) of Lemma \ref{L:good places} holds for some $v \in S_{good}\cap S_{\K_2}$.
\item
If $\KK_1/\F$ is not cyclic, then $\KK_1 = \K_3\K_4$ for some CM-extensions $\K_3/\F$ and $\K_4/\F$ inside $\KK_1$.
In this case, 
the triple product $L$-function factorizes into product of Rankin--Selberg $L$-functions for $\GL_2(\A_\F) \times \GL_2(\A_\F)$:
\[
L^S(s+\tfrac{1}{2},\itPi,{\rm As}) = L^S(s,\itPi_3' \times \itPi_2)L^S(s,\itPi_4' \times \itPi_2)
\]
for sufficiently large finite set $S$ of places, where $\itPi_i' = I_{\KK_i}^\F\left(\chi_1\vert_{\A_{\KK_i}^\times}\right) \otimes |\mbox{ }|_{\A_\F}^{-1/2}$ for $i=3,4$.
Let $\eta$ be a finite order Hecke character of $\A_\E^\times$. Then similarly we have
\[
L^S(s+\tfrac{1}{2},\itPi\otimes\eta,{\rm As}) = L^S(s,\itPi_3' \times \itPi_2 \otimes \eta\vert_{\A_\F^\times})L^S(s,\itPi_4' \times \itPi_2 \otimes \eta\vert_{\A_\F^\times}).
\]
Note that $L(1-\tfrac{{\sf w}}{2},\itPi,{\rm As})\neq 0$. Indeed, the Rankin--Selberg $L$-functions are non-vanishing at $s=\tfrac{1-{\sf w}}{2}$ by the result of Shahidi \cite[Theorem 5.2]{Shahidi1981}.
Note that $L(s,\itPi_i' \times \itPi_2 \otimes \eta\vert_{\A_\F^\times})$ does not necessary factorize into Hecke $L$-functions, unless $\K_i=\K_2$. Instead of using Lemma \ref{L:period relation}, we invoke the result of Shimura \cite[Theorem 4.2]{Shimura1978} on the algebraicity of critical values of Rankin--Selberg $L$-functions.
Note that $\itPi_3'$ and $\itPi_4'$ are cohomological irreducible cuspidal automorphic representations of $\GL_2(\A_\F)$.
By the theorem of Shimura, we have
\[
\sigma \left( \frac{L^{(\infty)}(\tfrac{1-{\sf w}}{2},\itPi_i' \times \itPi_2 \otimes \eta\vert_{\A_\F^\times})}{G(\eta\vert_{\A_\F^\times})^2\cdot L^{(\infty)}(\tfrac{1-{\sf w}}{2},\itPi_i' \times \itPi_2)}\right)  = \frac{L^{(\infty)}(\tfrac{1-{\sf w}}{2},{}^\sigma\!\itPi_i' \times {}^\sigma\!\itPi_2 \otimes {}^\sigma\!\eta\vert_{\A_\F^\times})}{G({}^\sigma\!\eta\vert_{\A_\F^\times})^2\cdot L^{(\infty)}(\tfrac{1-{\sf w}}{2},{}^\sigma\!\itPi_i' \times {}^\sigma\!\itPi_2)}
\]
for all $\sigma\in{\rm Aut}(\C)$ and $i=3,4$.
Here we have used the balanced condition (\ref{E:balanced 2}) to guarantee that $\tfrac{1-{\sf w}}{2}$ is a critical point for the Rankin--Selberg $L$-functions.
Therefore, (\ref{E:Case 1 proof 2}) also holds. 
By (\ref{E:Case 1 proof 1}) and (\ref{E:Case 1 proof 2}), the algebraicity holds for $L^{(\infty)}(1-\tfrac{\sf w}{2},\itPi,{\rm As})$ if and only if it holds for $L^{(\infty)}(1-\tfrac{\sf w}{2},\itPi\otimes\eta,{\rm As})$.
Assume further $\eta$ is chosen so that $\eta\vert_{\A_\F^\times}$ has prime order and $[\L:\F]$ is coprime to $[\K_1\K_2\L_\omega:\F]$. Here $\L/\F$ is the cyclic Galois extension associated $\eta\vert_{\A_\F^\times}$ by class field theory.
Then $\K_1\K_2\L_\omega \cap \L = \F$ and the following set of places of $\F$ has positive density by CDT for $\K_1\K_2\L_\omega\L/\F$:
\[
\left\{v \nmid 2\cdot\infty \,\left\vert\, \mbox{$v$ splits in $\K_1\K_2\L_\omega$, $\itPi_v$ and $\eta_v$ are unramified, and $\eta\vert_{\F_v^\times} \neq 1$}\right\}\right..
\]
It then follows from Lemma \ref{L:ST} for $\itPi_2$ that there exists $v$ belongs to the above set such that ${\rm tr}(C(\itPi_{2,v})) \neq 0$. On the other hands, it is clear that the above set is contained in the set (\ref{E:good places}) with $\itPi$ and $\omega$ replaced by $\itPi\otimes \eta$ and $\omega\eta\vert_{\A_\F^\times}^2$. 
Indeed, since $v$ splits in $\L_\omega$ and $\eta\vert_{\F_v^\times}$ has odd order, we have $\omega_v\eta\vert_{\F_v^\times}^2 = |\mbox{ }|^{\sf w}_v\eta\vert_{\F_v^\times}^2 \neq |\mbox{ }|_v^{\sf w}$.
Thus (1) of Lemma \ref{L:good places} holds with $\itPi$ replaced by $\itPi\otimes\eta$, and the algebraicity holds for $L^{(\infty)}(1-\tfrac{\sf w}{2},\itPi\otimes\eta,{\rm As})$.
\end{itemize}

\subsubsection{\underline{Case $3$}}\label{SS:Case 3}

We consider two subcases:
\begin{itemize}
\item[\underline{3-(i)}] $\itPi$ is non-CM.
\item[\underline{3-(ii)}] There exists CM-extension $\KK'/\E$  and algebraic Hecke character $\chi$ of $\A_{\KK'}^\times$ such that $\itPi = I_{\KK'}^\E(\chi) \otimes |\mbox{ }|_{\A_\E}^{-1/2}$.
\end{itemize}
The proof for Case is similar to that of Case $1$-(i) by applying Lemmas \ref{L:density} and \ref{L:ST}.

In Case 3-(ii), let $\widetilde{\K}'$ be the Galois closure of $\K'/\F$. We claim that there exists a quadratic extension $\K/\F$ inside $\K'$.
The assertion is clear if $\E/\F$ is Galois. When $\E/\F$ is not Galois, $\widetilde{\K}' = \widetilde{\E}\K'$ has degree $12$ over $\F$, where $\widetilde{\E}$ is the Galois closure of $\E/\F$. Thus $\Gal(\widetilde{\K}'/\F)$ is isomorphic to either the Dihedral group $D_6$ or the alternating group $A_4$ or the generalized quaternion group $Q_{12}$. Since $\E/\F$ is not Galois, there is a quadratic extension over $\F$ inside $\widetilde{\E}$. 
On the other hand, $A_4$ has no subgroup of index $2$.
Also $Q_{12}$ has a unique element of order $2$ whereas $\widetilde{\E}/\F$ and $\K'/\F$ are distinct extensions inside $\widetilde{\K}'$ of degree $6$. 
Thus $\Gal(\widetilde{\K}'/\F)$ is isomorphic to neither $A_4$ nor $Q_{12}$.
In $D_6$, each index $6$ subgroup is contained in a subgroup of index $2$. This establishes the claim.
Note that $\K/\F$ must be a CM-extension, since $\K'/\F$ it totally imaginary.
Then $\KK' = \E\K$ and the triple product $L$-function factorizes into Hecke $L$-functions (cf.\,\cite[(2.5)]{Ikeda1992}):
\[
L^S(s+\tfrac{1}{2},\itPi,{\rm As}) = L^S(s-1,\,\chi\vert_{\A_\KK^\times})L^S(s-1,\,(\chi\circ{\rm N}_{\KK'/\KK})\chi^{-1}\chi^c)
\]
for sufficiently large finite set $S$ of places. Here $c$ is the non-trivial automorphism in $\Gal(\KK'/\E)$.
Let $\eta$ be a finite order Hecke character of $\A_\F^\times$ with parallel signature $1$. Then $\itPi\otimes(\eta\circ{\rm N}_{\E/\F}) = I_{\K'}^\E(\chi\cdot(\eta\circ{\rm N}_{\K'/\F})) \otimes |\mbox{ }|_{\A_\E}^{-1/2}$ and similarly we have
\[
L^S(s+\tfrac{1}{2},\itPi\otimes(\eta\circ{\rm N}_{\E/\F}),{\rm As}) = L^S(s-1,\,\chi\vert_{\A_\KK^\times}\cdot(\eta\circ{\rm N}_{\K/\F})^3)L^S(s-1,\,(\chi\circ{\rm N}_{\KK'/\KK})\chi^{-1}\chi^c\cdot(\eta\circ{\rm N}_{\K'/\F})^3).
\]
The unbalanced condition (\ref{E:balanced 2}) implies that the characters $\chi\vert_{\A_\KK^\times}$ and $(\chi\circ{\rm N}_{\KK'/\KK})\chi^{-1}\chi^c$ satisfy the assumption in Lemma \ref{L:period relation} for the CM-extensions $\K/\F$ and $\K'/\E$, respectively.
Thus $L(1-\tfrac{\sf w}{2},\itPi,{\rm As}) \neq 0$ and we have
\begin{align}\label{E:Case 1 proof 3}
\sigma \left( \frac{L^{(\infty)}(1-\tfrac{\sf w}{2},\itPi\otimes(\eta\circ{\rm N}_{\K'/\F}),{\rm As})}{G(\eta)^3\cdot G(\eta\circ{\rm N}_{\E/\F})^3\cdot L^{(\infty)}(1-\tfrac{\sf w}{2},\itPi,{\rm As})}\right) = \frac{L^{(\infty)}(1-\tfrac{\sf w}{2},{}^\sigma\!\itPi\otimes({}^\sigma\!\eta\circ{\rm N}_{\K'/\F}),{\rm As})}{G({}^\sigma\!\eta)^3\cdot G({}^\sigma\!\eta\circ{\rm N}_{\E/\F})^3\cdot L^{(\infty)}(1-\tfrac{\sf w}{2},{}^\sigma\!\itPi,{\rm As})}
\end{align}
for all $\sigma \in {\rm Aut}(\C)$.
It follows easily from the property (\ref{E:Galois Gauss sum}) of Gauss sum that
\begin{align}\label{E:Case 1 proof 4}
\sigma \left(  \frac{G(\eta\circ{\rm N}_{\E/\F})}{G(\eta)^3}\right) = \frac{G({}^\sigma\!\eta\circ{\rm N}_{\E/\F})}{G({}^\sigma\!\eta)^3}
\end{align}
for all $\sigma \in {\rm Aut}(\C)$.
By (\ref{E:Case 1 proof 1}), (\ref{E:Case 1 proof 3}), and (\ref{E:Case 1 proof 4}), the algebraicity holds for $L^{(\infty)}(1-\tfrac{\sf w}{2},\itPi,{\rm As})$ if and only if it holds for $L^{(\infty)}(1-\tfrac{\sf w}{2},\itPi\otimes(\eta\circ{\rm N}_{\K'/\F}),{\rm As})$.
Assume $\eta$ has prime order and $\L/\F$ be the totally real Galois extension associated to $\eta$ by class field theory. 
We assume further that $\eta$ is chosen so that
$[\L:\F]$ is coprime to $[\widetilde{\K}'\L_\omega:\F]$.
In particular, $\widetilde{\K}'\L_\omega\cap \L = \F$.
Therefore, by CDT for $\widetilde{\K}'\L_\omega\L/\F$, the following set of places of $\F$ has positive density:
\[
\left\{ v \nmid2\cdot\infty\,\left\vert\,\mbox{$v$ splits in $\widetilde{\K}'\L_\omega$, $\itPi_v$ and $\eta_v$ are unramified, and $\eta_v \neq 1$}\right\}\right..
\]
It then follows from Lemma \ref{L:ST} for $\itPi$ that there exists $v$ belongs to the above set such that ${\rm tr}(C(\itPi_{v}^{(i)})) \neq 0$ for some $1 \leq i \leq 3$. On the other hands, it is clear that the above set is contained in the set (\ref{E:good places}) with $\itPi$ and $\omega$ replaced by $\itPi\otimes (\eta\circ{\rm N}_{\K'/\F})$ and $\omega\eta^6$. 
Indeed, since $v$ splits in $\L_\omega$ and $[\L:\F]$ is coprime to $6$, we have $\omega_v\eta_v^6 = |\mbox{ }|^{\sf w}_v\eta_v^6 \neq |\mbox{ }|_v^{\sf w}$.
Thus (1) of Lemma \ref{L:good places} holds with $\itPi$ replaced by $\itPi\otimes(\eta\circ{\rm N}_{\K'/\F})$, and the algebraicity holds for $L^{(\infty)}(1-\tfrac{\sf w}{2},\itPi\otimes(\eta\circ{\rm N}_{\K'/\F}),{\rm As})$.
This completes the proof of Theorem \ref{T:main 3} for $m+\tfrac{1}{2} \in {\rm Crit}^+(\itPi,{\rm As})$.

\subsubsection{Left-half critical values}\label{SS:LHC}

In this section we prove the algebraicity for left-half critical $L$-values.  Let $m+\tfrac{1}{2}$ be a left-half critical point for $L(s,\itPi,{\rm As})$. Then $-m+\tfrac{1}{2}$ is a right-half critical point for $L(s,\itPi^\vee,{\rm As})$. By the functional equation (\ref{E:fe}), we have
\[
L^{(\infty)}(m+\tfrac{1}{2},\itPi,{\rm As}) = \prod_{v \nmid \infty}\varepsilon(m+\tfrac{1}{2},\itPi_v,{\rm As},\psi_v)\prod_{v \in S_\infty}\gamma(m+\tfrac{1}{2},\itPi_v,{\rm As},\psi_v)\cdot L^{(\infty)}(-m+\tfrac{1}{2},\itPi^\vee,{\rm As}).
\]
Here we take $\psi_\F = \bigotimes_v\psi_v$ be the standard additive character of $\F$. 
By (\ref{E:Galois Gauss sum}) , (\ref{E:local factor archimedean}), and (\ref{E:local factor non-archimedean}), we have
\begin{align*}
\prod_{v \in S_\infty}\gamma(m+\tfrac{1}{2},\itPi_v,{\rm As},\psi_v) &\in (2\pi\sqrt{-1})^{8dm+4d{\sf w}}\cdot\Q^\times,\\
\sigma \left( \frac{\prod_{v \nmid \infty}\varepsilon(m+\tfrac{1}{2},\itPi_v,{\rm As},\psi_v)}{G(\omega)^4}\right) &= \frac{\prod_{v \nmid \infty}\varepsilon(m+\tfrac{1}{2},{}^\sigma\!\itPi_v,{\rm As},\psi_v)}{G({}^\sigma\!\omega)^4}
\end{align*}
for all $\sigma\in{\rm Aut}(\C)$. 
Also note that
\[
p(\itPi,{\rm As}) \in (2\pi\sqrt{-1})^{4d{\sf w}}\cdot G(\omega)^4\cdot p(\itPi^\vee,{\rm As}) \cdot \Q^\times
\]
Therefore, the algebraicity of $L^{(\infty)}(m+\tfrac{1}{2},\itPi,{\rm As})$ follows from that of $L^{(\infty)}(-m+\tfrac{1}{2},\itPi^\vee,{\rm As})$.
This completes the proof of Theorem \ref{T:main 3}.

\subsection{Proof of Theorem \ref{T:DC}}

We keep the notation of \S\,\ref{SS:symmetric cube}.
In this section, we prove Theorem \ref{T:DC}.
First we recall a result of Shimura \cite[Theorem 4.3]{Shimura1978} on the algebraicity of twisted standard $L$-function of $\itPi$.
\begin{thm}[Shimura]\label{T:Shimura}
There exists a sequence of non-zero complex numbers $(p({}^\sigma\!\itPi,\underline{\varepsilon}))_{\sigma \in {\rm Aut}(\C)}$ defined for each $\underline{\varepsilon} \in \{\pm1\}^{S_\infty}$ satisfying the following properties:
\begin{itemize}
\item[(i)] Let $\chi$ be a finite order Hecke character of $\A_\F^\times$ and $m+\tfrac{1}{2} \in \Z+\tfrac{1}{2}$ be a critical point of the tiwsted standard $L$-function $L(s,\itPi \otimes \chi)$. We have
\begin{align*}
&\sigma\left(\frac{L^{(\infty)}(m+\tfrac{1}{2},\itPi \otimes \chi)}{(2\pi\sqrt{-1})^{dm+\sum_{v \in S_\infty}(\kappa_v+{\sf w})/2}\cdot G(\chi) \cdot p(\itPi,(-1)^m\cdot{\rm sgn}(\chi))}\right)\\
& = \frac{L^{(\infty)}(m+\tfrac{1}{2},{}^\sigma\!\itPi \otimes {}^\sigma\!\chi)}{(2\pi\sqrt{-1})^{dm+\sum_{v \in S_\infty}(\kappa_v+{\sf w})/2}\cdot G({}^\sigma\!\chi) \cdot p({}^\sigma\!\itPi,(-1)^m\cdot{\rm sgn}(\chi))}
\end{align*}
for all $\sigma \in {\rm Aut}(\C)$.
\item[(ii)] Let $\underline{\varepsilon}\in \{\pm1\}^{S_\infty}$.
We have
\begin{align*}
\sigma \left(\frac{p(\itPi,\underline{\varepsilon})\cdot p(\itPi,-\underline{\varepsilon})}{(2\pi\sqrt{-1})^d\cdot(\sqrt{-1})^{d{\sf w}}\cdot G(\omega_\itPi)\cdot\Vert f_\itPi \Vert}\right) = \frac{p({}^\sigma\!\itPi,\underline{\varepsilon})\cdot p({}^\sigma\!\itPi,-\underline{\varepsilon})}{(2\pi\sqrt{-1})^d\cdot(\sqrt{-1})^{d{\sf w}}\cdot G({}^\sigma\!\omega_\itPi)\cdot\Vert f_{{}^\sigma\!\itPi} \Vert}
\end{align*}
for all $\sigma \in {\rm Aut}(\C)$.
\end{itemize}

\end{thm}

\begin{rmk}\label{R:Shimura}
Let ${\bf f}$ be the Hilbert modular cusp newform associated to $\itPi$. In the notation of \cite[Theorem 4.3]{Shimura1978}, we have
\[
p(\itPi,\underline{\varepsilon}) = (2\pi\sqrt{-1})^{\sum_{v \in S_\infty}(\kappa_v-\kappa_0)/2}\cdot u((-1)^{(\kappa_0+{\sf w})/2}\underline{\varepsilon},{\bf f}),
\]
where $\kappa_0 = \max_{v \in S_\infty}\{\kappa_v\}$.
\end{rmk}

We begin the proof of Theorem \ref{T:DC}. Let $m+\tfrac{1}{2}$ be a critical point for $L(s,\itPi,{\rm Sym}^3\otimes\chi)$ with $m+\tfrac{3{\sf w}}{2} \neq 0$. Note that the existence of such critical point is equivalent to our assumption that $\kappa_v \geq 3$ for all $v \in S_\infty$. 
For irreducible cuspidal automorphic representation $\itPi_1\boxtimes\itPi_2 \boxtimes \itPi_3$ on $\GL_2(\A_\E)$ with $\E=\F\times\F\times\F$, it is customary to write 
\[
L(s,\itPi_1 \times \itPi_2 \times \itPi_3) = L(s,\itPi_1\boxtimes\itPi_2 \boxtimes \itPi_3,{\rm As})
\]
for the triple product $L$-function. 
We take $(\itPi_1,\itPi_2,\itPi_3) = (\itPi,\itPi,\itPi\otimes\chi)$. In this case, we have the factorization of $L$-functions:
\[
L(s,\itPi \times \itPi \times \itPi\otimes\chi) = L(s,\itPi,{\rm Sym}^3\otimes\chi)L(s,\itPi \otimes \chi\omega_\itPi)^2.
\]
Note that the assumption $m+\tfrac{3{\sf w}}{2} \neq 0$ implies that $L(m+\tfrac{1}{2},\itPi \otimes \chi\omega_\itPi) \neq 0$. 
Indeed, for $m+\tfrac{3{\sf w}}{2} = \pm\tfrac{1}{2}$, the non-vanishing follows from \cite{JS1976}.
For $m+\tfrac{3{\sf w}}{2}\geq 1$, the Euler product defining $L(s,\itPi\otimes \chi\omega_\itPi)$ converges absolutely at $s=m+\tfrac{1}{2}$.
For $m+\tfrac{3{\sf w}}{2}\leq -1$, the non-vanishing then follows from the functional equation.
Therefore, we easily deduce Conjecture \ref{C:DC} for $L^{(\infty)}(m+\tfrac{1}{2},\itPi,{\rm Sym}^3\otimes\chi)$ from Theorems \ref{T:main 3} and \ref{T:Shimura}.
Now we consider the central critical point $s=\tfrac{1-3{\sf w}}{2}$. Note that in this case ${\sf w}$ is even. 
We may assume $\itPi$ is non-CM.
Similarly as above, we have $L(\tfrac{3-3{\sf w}}{2},\itPi,{\rm Sym}^3\otimes\chi) \neq 0$.
Let $\itSigma$ be the functorial lift of $\itPi$ to $\GL_4(\A_\F)$ with respect to the symmetric cube representation (cf.\,\cite{KS2002}). 
Then we have
\[
L(s,\itSigma\otimes\chi) = L(s,\itPi,{\rm Sym}^3\otimes\chi).
\]
Note that $\itSigma$ is a cohomological irreducible cuspidal automorphic representation of $\GL_4(\A_\F)$. We have the factorization of the twisted exterior square $L$-function
\[
L(s,\itSigma,\wedge^2 \otimes \omega_\itPi^{-3}) = L(s,\itPi,{\rm Sym}^4\otimes \omega_\itPi^{-2})\cdot \zeta_\F(s),
\]
where $L(s,\itPi,{\rm Sym}^4\otimes \omega_\itPi^{-2})$ is the twisted symmetric fourth $L$-function of $\itPi$ by $\omega_\itPi^{-2}$.
In particular, $L(s,\itSigma,\wedge^2 \otimes \omega_\itPi^{-3})$ has a pole at $s=1$ as $L(s,\itPi,{\rm Sym}^4\otimes \omega_\itPi^{-2})$ does not vanish at $s=1$.
In this case, for $\sigma \in {\rm Aut}(\C)$ and $\underline{\varepsilon} \in \{\pm1\}^{S_\infty}$, we have the Betti--Shilika period $\omega^{\underline{\varepsilon}}({}^\sigma\!\itSigma_f) \in \C^\times$ in \cite[Definition/Proposition 4.2.1]{GR2014}. 
Fix a finite order Hecke character $\eta$ of $\A_\F^\times$ such that ${\rm sgn}(\eta) = - {\rm sgn}(\chi)$.
By \cite[Theorem 7.1.2]{GR2014}, we have
\begin{align}\label{E:main 2 proof 1}
\begin{split}
&\sigma \left(\frac{ L^{(\infty)}(m+\tfrac{1}{2},\itSigma \otimes \chi)L^{(\infty)}(m+\tfrac{1}{2},\itSigma \otimes \eta)^{-1}}{ \omega^{(-1)^m{\rm sgn}(\chi)}(\itSigma_f)\omega^{(-1)^{m+1}{\rm sgn}(\chi)}(\itSigma_f)^{-1}\cdot G(\chi\eta^{-1})^2}\right)\\
& = \frac{L^{(\infty)}(m+\tfrac{1}{2},{}^\sigma\!\itSigma \otimes {}^\sigma\!\chi)L^{(\infty)}(m+\tfrac{1}{2},{}^\sigma\!\itSigma \otimes {}^\sigma\!\eta)^{-1}}{\omega^{(-1)^m{\rm sgn}(\chi)}({}^\sigma\!\itSigma_f)\omega^{(-1)^{m+1}{\rm sgn}(\chi)}({}^\sigma\!\itSigma_f)^{-1}\cdot G({}^\sigma\!\chi{}^\sigma\!\eta^{-1})^2}
\end{split}
\end{align}
for all critical points $m+\tfrac{1}{2} \neq \tfrac{1-3{\sf w}}{2}$ and $\sigma \in {\rm Aut}(\C)$. 
On the other hand, by \cite[Theorem 7.40]{HR2020}, we have
\begin{align}\label{E:main 2 proof 2}
\begin{split}
&\sigma \left(\frac{L^{(\infty)}(m-\tfrac{1}{2},\itSigma \otimes \chi)}{(2\pi\sqrt{-1})^{-2d}\cdot\Omega^{(-1)^m{\rm sgn}(\chi)}(\itSigma_f)\cdot L^{(\infty)}(m+\tfrac{1}{2},\itSigma \otimes \chi)} \right)\\
& = \frac{L^{(\infty)}(m-\tfrac{1}{2},{}^\sigma\!\itSigma \otimes {}^\sigma\!\chi)}{(2\pi\sqrt{-1})^{-2d}\cdot\Omega^{(-1)^m{\rm sgn}(\chi)}({}^\sigma\!\itSigma_f)\cdot L^{(\infty)}(m+\tfrac{1}{2},{}^\sigma\!\itSigma \otimes {}^\sigma\!\chi)}
\end{split}
\end{align}
for all critical points $m+\tfrac{1}{2} \neq \tfrac{1-3{\sf w}}{2}$ and $\sigma \in {\rm Aut}(\C)$, 
where $\Omega^{\underline{\varepsilon}}(\itSigma_f) \in \C^\times$ is the relative period defined in \cite[Definition 5.10]{HR2020}. As explained in \cite[\S\,5.2.3]{HR2020}, we have
\[
\sigma \left( \Omega^{\underline{\varepsilon}}(\itSigma_f)\cdot \frac{\omega^{\underline{\varepsilon}}(\itSigma_f)}{\omega^{-\underline{\varepsilon}}(\itSigma_f)} \right) = \Omega^{\underline{\varepsilon}}({}^\sigma\!\itSigma_f)\cdot \frac{\omega^{\underline{\varepsilon}}({}^\sigma\!\itSigma_f)}{\omega^{-\underline{\varepsilon}}({}^\sigma\!\itSigma_f)}
\]
for all $\sigma \in {\rm Aut}(\C)$.
We thus conclude from (\ref{E:main 2 proof 1}) and (\ref{E:main 2 proof 2}) that
\begin{align*}
\sigma \left( \frac{L^{(\infty)}(\tfrac{1-3{\sf w}}{2},\itSigma \otimes \chi)}{(2\pi\sqrt{-1})^{-2d}\cdot G(\chi\eta^{-1})^2\cdot L^{(\infty)}(\tfrac{3-3{\sf w}}{2},\itSigma \otimes \eta)}\right) = \frac{L^{(\infty)}(\tfrac{1-3{\sf w}}{2},{}^\sigma\!\itSigma \otimes {}^\sigma\!\chi)}{(2\pi\sqrt{-1})^{-2d}\cdot G({}^\sigma\!\chi{}^\sigma\!\eta^{-1})^2\cdot L^{(\infty)}(\tfrac{3-3{\sf w}}{2},{}^\sigma\!\itSigma \otimes {}^\sigma\!\eta)}
\end{align*}
for all $\sigma \in {\rm Aut}(\C)$.
Therefore, Conjecture \ref{C:DC} for $L^{(\infty)}(\tfrac{1-3{\sf w}}{2},\itSigma \otimes \chi)$ follows from that of $L^{(\infty)}(\tfrac{3-3{\sf w}}{2},\itSigma \otimes \eta)$.
This completes the proof.

\section{Computation of local zeta integrals}\label{S:local zeta integral}

In this section, we compute the non-archimedean local zeta integrals. The main result is Lemma \ref{L:local zeta 3}.
We follow the notation in (\ref{E:group elements}).

Let $\F$ be a non-archimedean local field of characteristic zero. Let $\frak{o}$, $\varpi$, and $q$ be the maximal compact subring of $\F$, a generator of the maximal ideal of $\frak{o}$, and the cardinality of $\frak{o} / \varpi\frak{o}$, respectively. Let $|\mbox{ }|$ be the absolute value on $\F$ normalized so that $|\varpi| = q^{-1}$. Fix a non-trivial additive character $\psi$ of $\F$ with conductor $\o$. 
Let $\mathcal{S}(\F)$ be the space of Schwartz functions on $\F$ which consisting of locally constant functions on $\F$ with compact support.
For $\varphi\in\mathcal{S}(\F)$, let $\widehat{\varphi}$ be the Fourier transform with respect to $\psi$ defined by
\[
\widehat{\varphi}(x) = \int_{\F}\varphi(y)\psi(xy)\,dy.
\]

Let $\itPi_1$, $\itPi_2$, and $\itPi_3$ be irreducible generic unramified representations of $\GL_2(\F)$ with central characters $\omega_{\itPi_1}$,  $\omega_{\itPi_2}$, and $\omega_{\itPi_3}$. Put $\omega = \omega_{\itPi_1}\omega_{\itPi_2}\omega_{\itPi_3}$.
For $i=1,2,3$, let $\alpha_{i},\beta_i$ be the Satake parameters of $\itPi_i$ and define unramified characters $\chi_i$ and $\mu_i$ of $\F^\times$ so that
\[
\chi_i(\varpi) = \alpha_i,\quad \mu_i(\varpi) = \beta_i.
\]
For Whittaker function $W_i \in \mathcal{W}(\itPi_i,\psi)$ with $i=1,2,3$ and meromorphic section $f^{(s)}$ of $I(\omega,s)$, recall the local zeta integral $Z(W_1\otimes W_2 \otimes W_3,f^{(s)})$ defined in (\ref{E:local zeta}) with $\eta_\F$ replace by
\[
\eta = \bp 0&0&0&-1&0&0 \\ 0&1&0&0&0&0 \\ 0&0&1&0&0&0 \\ 1&1&1&0&0&0 \\ 0&0&0&-1&1&0 \\ 0&0&0&-1&0&1\ep.
\]
Here we identify ${\bf G}(\F) = \{(g_1,g_2,g_3) \in \GL_2(\F)^3\,\vert\,\det g_1 = \det g_2 = \det g_3\}$ with subgroup of $\GSp_6(\F)$ by the embedding (\ref{E:embedding}).

For $\itPhi \in \mathcal{S}({\rm Sym}_3(\F))$, recall $f_\itPhi^{(s)}$ is the rational section of $I(\omega,s)$ such that $f_\itPhi^{(s)}$ is supported in $P_3(\F)J_3P_3(\F)$ and
\[
f_\itPhi^{(s)}(J_3{\bf n}(x)) = \itPhi(x).
\]

\begin{lemma}\label{L:local zeta 1}
Let $\itPhi \in \mathcal{S}({\rm Sym}_3(\F))$ such that $\itPhi$ factors through the quotient ${\rm Sym}_3(\varpi^{-n}\o) / {\rm Sym}_3(\varpi^n\o)$ for some $n > 0$.
Then $f_\itPhi^{(s)}$ is right invariant by ${\bf n}^-({\rm Sym}_3(\varpi^{2n}\o))$.
\end{lemma}

\begin{proof}
Since the support of $\itPhi$ is contained in ${\rm Sym}_3(\varpi^{-n}\o)$, it suffices to show that
\[
f_\itPhi^{(s)}(J_3{\bf n}(x){\bf n}^-(y)) = f_\itPhi^{(s)}(J_3{\bf n}(x))
\]
for all $x \in {\rm Sym}_3(\varpi^{-n}\o)$ and $y \in {\rm Sym}_3(\varpi^{2n}\o)$.
Note that 
\[
J_3{\bf n}(x){\bf n}^-(y) = \bp ({\bf 1}_3+yx)^{-1} & -y \\ 0 & {\bf 1}_3+xy\ep J_3{\bf n}(({\bf 1}_3+xy)^{-1}x).
\]
Since $\det({\bf 1}_3+xy) \in \o^\times$, we have
\[
f_\itPhi^{(s)}(J_3{\bf n}(x){\bf n}^-(y)) = f_\itPhi^{(s)}(J_3{\bf n}(({\bf 1}_3+xy)^{-1}x)).
\]
As $({\bf 1}_3+xy)^{-1}x - x \in {\rm Sym}_3(\varpi^n \o)$, we conclude from the condition on $\itPhi$ that
\[
f_\itPhi^{(s)}(J_3{\bf n}(({\bf 1}_3+xy)^{-1}x)) = f_\itPhi^{(s)}(J_3{\bf n}(x)).
\]
This completes the proof.
\end{proof}
For $n \geq 1$ and $1 \leq i \leq 3$, let $f_i^{(n)}$ be the section in the induced representation ${\rm Ind}_{P_1(\F)}^{\GL_2(\F)}(\chi_i\boxtimes\mu_i)$ such that the support of $f_i^{(n)}$
 is contained in $P_1(\F)J_1P_1(\F)J_1$ and
\[
 f_i^{(n)}({\bf n}^-(x)) = \mathbb{I}_{\varpi^{2n}\o}(x).
\]
Let $W_i^{(n)}\in \mathcal{W}(\itPi_i,\psi)$ be the Whittaker function of $\itPi_i$ defined by
\begin{align}\label{E:Whittaker function}
W_i^{(n)}(g) = \int_{\F}f_i^{(n)}(J_1{\bf n}(x)g)\overline{\psi(x)}\,dx.
\end{align}
Then it is easy to verify that
\[
W_i^{(n)}({\bf t}(a)J_1) = q^{-2n}\mu_i(a)|a|^{1/2}\mathbb{I}_{\varpi^{-2n}\o}(a).
\]

\begin{lemma}\label{L:local zeta 2}
For $W_1 \in \mathcal{W}(\itPi_1,\psi)$ and $\itPhi \in \mathcal{S}({\rm Sym}_3(\F))$ which factors through ${\rm Sym}_3(\varpi^{-n}\o) / {\rm Sym}_3(\varpi^n\o)$ for some $n >0$ with
\[
\itPhi\left(\bp x_{11} & x_{12} & x_{13} \\ x_{12} & x_{22} & x_{23} \\ x_{13} & x_{23} & x_{33} \ep\right) = \prod_{1 \leq i \leq j \leq 3}\varphi_{ij}(x_{ij}),
\]
we have
\begin{align*}
Z(W_1\otimes W_2^{(n)}\otimes W_3^{(n)},f_\itPhi^{(s)}) 
& = q^{-4n}(1+q^{-1})^{-2} \widehat{\varphi}_{22}(0)\widehat{\varphi}_{33}(0) \gamma(s+\tfrac{1}{2},\itPi_1\otimes\mu_2\mu_3,\psi)^{-1}\\
&\times\int_{\F}dx\int_{(\F^\times)^3}d^\times a_1\,d^\times a_2\,d^\times a_3\,W_1({\bf t}(a_1a_2a_3)J_1{\bf n}(x))\omega_{\itPi_1}^{-1}\mu_2^{-1}\mu_3^{-1}(a_1)|a_1|^{-s}\widehat{\varphi}_{23}(a_1)\\
&\quad\quad\quad\quad\quad\quad\quad\quad\quad\quad\quad\quad\quad\times \mu_2\chi_3(a_2)|a_2|^{s}
\mu_3\chi_2(a_3)|a_3|^{s}\varphi_{12}(a_2)\varphi_{13}(a_3)\varphi_{11}(x).
\end{align*}
Here $\gamma(s,\itPi_1\otimes\mu_2\mu_3,\psi)$ is the $\gamma$-factor of $\itPi_1\otimes\mu_2\mu_3$ with respect to $\psi$.
\end{lemma}

\begin{proof}
Rewrite $W_2^{(n)}$ and $W_3^{(n)}$ in terms of Whittaker integral (\ref{E:Whittaker function}), we have
\begin{align*}
&Z(W_1\otimes W_2^{(n)}\otimes W_3^{(n)},f_\itPhi^{(s)})\\
&= \int_{\F^\times U^0(\F) \backslash {\bf G}(\F)}W_1(g_1)W_2^{(n)}(g_2)W_3^{(n)}(g_3)f_\itPhi^{(s)}(\eta\cdot\iota(g_1,g_2,g_3))\,dg\\
&= \int_{\F^\times \backslash{\bf G}(\F)}W_1(g_1)f_2^{(n)}(J_1g_2)f_3^{(n)}(J_1g_3)f_\itPhi^{(s)}(\eta\cdot\iota(g_1,g_2,g_3))\,dg\\
& = \int_{\F^\times}d^\times a_1 \int_{\F^\times \backslash \SL_2(\F)}dg_1 \int_{\SL_2(\F)^2}dg_2\,dg_3\, W_1({\bf t}(a_1)g_1)f_2^{(n)}(J_1{\bf t}(a_1)g_2)f_3^{(n)}(J_1{\bf t}(a_1)g_3)\\
& \quad\quad\quad\quad\quad\quad\quad\quad\quad\quad\quad\quad\quad\quad\quad\quad\quad\quad\quad\quad\quad\times f_\itPhi^{(s)}(\eta\cdot\iota({\bf t}(a_1)g_1,{\bf t}(a_1)g_2,{\bf t}(a_1)g_3))\\
& = \int_{\F^\times}d^\times a_1 \int_{\F^\times \backslash \SL_2(\F)}dg_1 \, W_1({\bf t}(a_1)g_1) \mu_2\mu_3(a_1)|a_1|^s\\
&\quad\quad\quad\quad\times\int_{\SL_2(\F)^2}dg_2\,dg_3\,f_2^{(n)}(J_1g_2)f_3^{(n)}(J_1g_3)f_\itPhi^{(s)}(\eta\cdot\iota(g_1,g_2,g_3)).
\end{align*}
Recall we have the integration formula
\[
\int_{\SL_2(\F)}f(g)\,dg = (1+q^{-1})^{-1}\int_{\F^\times}d^\times a\int_{\F^2}dx\,dy \,f({\bf m}(a){\bf n}^-(x){\bf n}(y))
\]
for any integrable function $f$ on $\SL_2(\F)$. 
Thus we have
\begin{align*}
&Z(W_1\otimes W_2^{(n)}\otimes W_3^{(n)},f_\itPhi^{(s)})\\
&= (1+q^{-1})^{-2}\int_{(\F^\times)^3}d^\times a_1\,d^\times a_2\,d^\times a_3 \int_{\F^6}dx_1\,dy_1\,dx_2\,dy_2\,dx_3\,dy_3 W_1({\bf t}(a_1){\bf n}^-(x_1){\bf n}(y_1))\mu_2\mu_3(a_1)|a_1|^s\\
&\quad\quad\quad\quad\quad\quad\quad\quad\quad\quad\quad\quad\quad\quad\quad\quad\times f_2^{(n)}(J_1 {\bf m}(a_2){\bf n}^-(x_2)J_1{\bf n}^-(y_2))f_3^{(n)}(J_1 {\bf m}(a_3){\bf n}^-(x_3)J_1{\bf n}^-(y_3))\\
&\quad\quad\quad\quad\quad\quad\quad\quad\quad\quad\quad\quad\times f_\itPhi^{(s)}(\eta\cdot\iota({\bf n}^-(x_1){\bf n}(y_1),{\bf m}(a_2){\bf n}^-(x_2)J_1{\bf n}^-(y_2),{\bf m}(a_3){\bf n}^-(x_3)J_1{\bf n}^-(y_3)))\\
&= q^{-4n}(1+q^{-1})^{-2}\int_{(\F^\times)^3}d^\times a_1\,d^\times a_2\,d^\times a_3 \int_{\F^4}dx_1\,dy_1\,dx_2\,dx_3\, W_1({\bf t}(a_1){\bf n}^-(x_1){\bf n}(y_1))\mu_2\mu_3(a_1)|a_1|^s\\
&\quad\quad\quad\quad\quad\times \chi_2^{-1}\mu_2(a_2)|a_2|^{-1}\chi_3^{-1}\mu_3(a_3)|a_3|^{-1} f_\itPhi^{(s)}(\eta\cdot\iota({\bf n}^-(x_1){\bf n}(y_1),{\bf m}(a_2){\bf n}^-(x_2)J_1,{\bf m}(a_3){\bf n}^-(x_3)J_1)).
\end{align*}
Here the second equality follows from Lemma \ref{L:local zeta 1} and the definition of $f_2^{(n)}$ and $f_3^{(n)}$.
Indeed, 
\[
\int_{\F}f_i^{(n)}({\bf n}^-(y_i))\,dy_i = q^{-2n}
\]
for $i=2,3$.
Note that
\begin{align*}
&\eta\cdot\iota({\bf n}^-(x_1){\bf n}(y_1),{\bf m}(a_2){\bf n}^-(x_2)J_1,{\bf m}(a_3){\bf n}^-(x_3)J_1)\\
& = {\bf m}\left(\bp -1 & a_2x_1 & a_3x_1 \\ 0 & a_2 & 0 \\ 0 & 0 & a_3 \ep\right){\bf n}\left( \bp -x_1 & 0 & 0 \\ 0 & 0 & 0 \\ 0 & 0 & 0 \ep\right)J_3{\bf n}\left( \bp y_1 & a_2 & a_3 \\ a_2 & -x_2-a_2^2x_1 & -a_2a_3x_1 \\ a_3 & -a_2a_3x_1 & -x_3-a_3^2x_1 \ep\right).
\end{align*}
Thus 
\begin{align*}
&f_\itPhi^{(s)} (\eta\cdot\iota({\bf n}^-(x_1){\bf n}(y_1),{\bf m}(a_2){\bf n}^-(x_2)J_1,{\bf m}(a_3){\bf n}^-(x_3)J_1))\\
& = \omega(a_2a_3)|a_2a_3|^{2s+2}\varphi_{11}(y_1)\varphi_{22}(-x_2-a_2^2x_1)\varphi_{33}(-x_3-a_3^2x_1)\varphi_{12}(a_2)\varphi_{13}(a_3)\varphi_{23}(-a_2a_3x_1).
\end{align*}
We conclude that the integrations in $x_2$ and $x_3$ are equal to $\widehat{\varphi}_{22}(0)$ and $\widehat{\varphi}_{33}(0)$ and we have
\begin{align*}
&Z(W_1\otimes W_2^{(n)}\otimes W_3^{(n)},f_\itPhi^{(s)})\\
& = q^{-4n}(1+q^{-1})^{-2} \widehat{\varphi}_{22}(0)\widehat{\varphi}_{33}(0)\\
&\times\int_{(\F^\times)^3}d^\times a_1\,d^\times a_2\,d^\times a_3 \int_{\F^2}dx_1\,dy_1\,W_1({\bf t}(a_1){\bf n}^-(x_1){\bf n}(y_1))\mu_2\mu_3(a_1)|a_1|^s\\
& \quad\times \omega\chi_2^{-1}\mu_2(a_2)|a_2|^{2s+1}
\omega\chi_3^{-1}\mu_3(a_3)|a_3|^{2s+1}\varphi_{12}(a_2)\varphi_{13}(a_3) \varphi_{11}(y_1)\varphi_{23}(-a_2a_3x_1).
\end{align*}
Recall we have the functional equation 
\[
\int_{\F^\times} W({\bf t}(a))\chi(a)|a|^s\,d^\times a = \gamma(s+\tfrac{1}{2},\itPi_1\otimes\chi,\psi)^{-1}\int_{\F^\times} W({\bf t}(a)J_1)\omega_{\itPi_1}^{-1}\chi^{-1}(a)|a|^{-s}\,d^\times a
\]
for all $W \in \mathcal{W}(\itPi_1,\psi)$ and character $\chi$ of $\F^\times$.
Thus we have
\begin{align*}
&\int_{\F}dx_1\int_{\F^\times}d^\times a_1\, W_1({\bf t}(a_1){\bf n}^-(x_1){\bf n}(y_1))\mu_2\mu_3(a_1)|a_1|^s\varphi_{23}(-a_2a_3x_1)\\
& = \gamma(s+\tfrac{1}{2},\itPi_1\otimes\mu_2\mu_3,\psi)^{-1}\int_{\F}dx_1\int_{\F^\times}d^\times a_1\, W_1({\bf t}(a_1)J_1{\bf n}^-(x_1){\bf n}(y_1))\omega_{\itPi_1}^{-1}\mu_2^{-1}\mu_3^{-1}(a_1)|a_1|^{-s}\varphi_{23}(-a_2a_3x_1)\\
& = |a_2a_3|^{-1} \gamma(s+\tfrac{1}{2},\itPi_1\otimes\mu_2\mu_3,\psi)^{-1} \int_{\F^\times}W_1({\bf t}(a_1)J_1{\bf n}(y_1))\omega_{\itPi_1}^{-1}\mu_2^{-1}\mu_3^{-1}(a_1)|a_1|^{-s}\widehat{\varphi}_{23}(a_1a_2^{-1}a_3^{-1})\,d^\times a_1\\
& = \omega_{\itPi_1}^{-1}\mu_2^{-1}\mu_3^{-1}(a_2a_3)|a_2a_3|^{-s-1} \gamma(s+\tfrac{1}{2},\itPi_1\otimes\mu_2\mu_3,\psi)^{-1}\\
&\times \int_{\F^\times}W_1({\bf t}(a_1a_2a_3)J_1{\bf n}(y_1))\omega_{\itPi_1}^{-1}\mu_2^{-1}\mu_3^{-1}(a_1)|a_1|^{-s}\widehat{\varphi}_{23}(a_1)\,d^\times a_1.
\end{align*}
This completes the proof.
\end{proof}

Now we specify the choices of $W_1$ and $\itPhi$ in the above lemma and obtain explicit values of the corresponding local zeta integrals.
Let $W_1^\circ \in \mathcal{W}(\itPi_1,\psi)$ be the $\GL_2(\o)$-invariant Whittaker function normalized so that $W_1(1)=1$. It is well-known that
\[
W_1^\circ({\bf t}(a)) = |a|^{1/2}\frac{\chi_1(a\varpi)-\mu_1(a\varpi)}{\chi_1(\varpi)-\mu_1(\varpi)}\cdot \mathbb{I}_\o(a).
\]

\begin{lemma}\label{L:local zeta 3}
\noindent
\begin{itemize}
\item[(1)]
Let $\itPhi \in \mathcal{S}({\rm Sym}_3(\F))$ defined by
\[
\itPhi\left(\bp x_{11} & x_{12} & x_{13} \\ x_{12} & x_{22} & x_{23} \\ x_{13} & x_{23} & x_{33} \ep\right) = \widehat{\mathbb{I}}_{\o^\times}(x_{11}){\mathbb{I}}_{\varpi^{-1}\o}(x_{22}){\mathbb{I}}_{\varpi^{-1}\o}(x_{33}){\mathbb{I}}_{\varpi^{-2}\o^\times}(x_{12}){\mathbb{I}}_{\varpi^{-1}\o}(x_{13})\widehat{\mathbb{I}}_{\o^\times}(x_{23}).
\]
We have
\begin{align*}
Z(W_1^\circ\otimes W_2^{(2)}\otimes W_3^{(2)},f_\itPhi^{(s)}) &= q^{-6}(1+q^{-1})^{-2} (\alpha_2\beta_3)^{2}(\beta_2\alpha_3)^{-2}\gamma(s+\tfrac{1}{2},\itPi_1\times\itPi_2\otimes\mu_3)^{-1}.
\end{align*}
\item[(2)] 
Let $\itPhi' \in \mathcal{S}({\rm Sym}_3(\F))$ defined by
\[
\itPhi'\left(\bp x_{11} & x_{12} & x_{13} \\ x_{12} & x_{22} & x_{23} \\ x_{13} & x_{23} & x_{33} \ep\right) = {\mathbb{I}}_{\varpi^{-1}\o}(x_{11}){\mathbb{I}}_{\varpi^{-1}\o}(x_{22}){\mathbb{I}}_{\varpi^{-1}\o}(x_{33})\widehat{\mathbb{I}}_{\o^\times}(x_{12})\widehat{\mathbb{I}}_{\o^\times}(x_{13})\widehat{\mathbb{I}}_{\o^\times}(x_{23}).
\]
We have
\begin{align*}
Z(W_1^{(1)}\otimes W_2^{(1)}\otimes W_3^{(1)},f_{\itPhi'}^{(s)})
& = q^{-3}(1+q^{-1})^{-2}\\
&\times\gamma(s+\tfrac{1}{2},\itPi_1\otimes\mu_2\mu_3,\psi)^{-1}\gamma(s+\tfrac{1}{2},\mu_1\mu_2\chi_3,\psi)^{-1}\gamma(s+\tfrac{1}{2},\mu_1\mu_3\chi_2,\psi)^{-1}.
\end{align*}
\end{itemize}
\end{lemma}

\begin{proof}
We begin with the first assertion.
Note that
\[
\widehat{\mathbb{I}}_{\o^\times} =(1-q^{-1}){\mathbb{I}}_{\o}-q^{-1}{\mathbb{I}}_{\varpi^{-1}\o^\times}.
\]
Thus it is clear that $\itPhi$ factors through ${\rm Sym}_3(\varpi^{-2}\o) / {\rm Sym}_3(\varpi^2\o)$.
By Lemma \ref{L:local zeta 2}, we have
\begin{align*}
&Z(W_1^\circ\otimes W_2^{(2)}\otimes W_3^{(2)},f_\itPhi^{(s)}) \\
& = q^{-6}(1+q^{-1})^{-2} \gamma(s+\tfrac{1}{2},\itPi_1\otimes\mu_2\mu_3,\psi)^{-1}\\
&\times\int_{\F}dx\int_{(\F^\times)^2}d^\times a_2\,d^\times a_3\,W_1^\circ({\bf t}(a_2a_3)J_1{\bf n}(x))\mu_2\chi_3(a_2)|a_2|^s\mu_3\chi_2(a_3)|a_3|^s{\mathbb{I}}_{\varpi^{-2}\o^\times}(a_2){\mathbb{I}}_{\varpi^{-1}\o}(a_3)\widehat{\mathbb{I}}_{\o^\times}(x).
\end{align*}
Note that $W_1^\circ({\bf t}(a))$ is supported in $\o \setminus \{0\}$ and 
\[
J_1{\bf n}(x) = {\bf n}(-x^{-1}){\bf m}(-x^{-1}){\bf n}^-(x^{-1})
\]
for $x \in \F^\times$.
In particular, for $x \in \varpi^{-1}\o$, we have
\[
W_1^\circ({\bf t}(a)J_1{\bf n}(x)) =\begin{cases}
W_1^\circ({\bf t}(a)) & \mbox{ if $x \in \o$},\\
(\alpha_1\beta_1)^{-1}\psi(-ax^{-1}) W_1^\circ({\bf t}(a\varpi^2))& \mbox{ if $x \in \varpi^{-1}\o^\times$}.
\end{cases}
\]
Therefore, 
\begin{align*}
&\int_{\F}dx\int_{(\F^\times)^2}d^\times a_2\,d^\times a_3\,W_1^\circ({\bf t}(a_2a_3)J_1{\bf n}(x))\mu_2\chi_3(a_2)|a_2|^s\mu_3\chi_2(a_3)|a_3|^s{\mathbb{I}}_{\varpi^{-2}\o^\times}(a_2){\mathbb{I}}_{\varpi^{-1}\o}(a_3)\widehat{\mathbb{I}}_{\o^\times}(x)\\
& = (1-q^{-1})\int_{\F}dx\int_{(\F^\times)^2}d^\times a_2\,d^\times a_3\,W_1^\circ({\bf t}(a_2a_3))\mu_2\chi_3(a_2)|a_2|^s\mu_3\chi_2(a_3)|a_3|^s{\mathbb{I}}_{\varpi^{-2}\o^\times}(a_2){\mathbb{I}}_{\varpi^{-1}\o}(a_3)\\
&-(\alpha_1\beta_1)^{-1}\int_{\F}dx\int_{(\F^\times)^2}d^\times a_2\,d^\times a_3\,W_1^\circ({\bf t}(a_2a_3\varpi^2))\mu_2\chi_3(a_2)|a_2|^s\mu_3\chi_2(a_3)|a_3|^s\\
&\quad\quad\quad\quad\quad\quad\quad\quad\quad\quad\quad\quad\quad\quad\quad\quad\quad\quad\times{\mathbb{I}}_{\varpi^{-2}\o^\times}(a_2){\mathbb{I}}_{\varpi^{-1}\o}(a_3)\widehat{\mathbb{I}}_{\o^\times}(a_2a_3\varpi)\\
& = (1-q^{-1})\cdot I_1 -(\alpha_1\beta_1)^{-1}\cdot I_2.
\end{align*}
To compute the integrals $I_1$ and $I_2$, note that
\begin{align*}
\mathbb{I}_\o(a_2a_3){\mathbb{I}}_{\varpi^{-2}\o^\times}(a_2){\mathbb{I}}_{\varpi^{-1}\o}(a_3) &= {\mathbb{I}}_{\varpi^{-2}\o^\times}(a_2){\mathbb{I}}_{\varpi^2\o}(a_3),\\
\mathbb{I}_\o(a_2a_3\varpi^2){\mathbb{I}}_{\varpi^{-2}\o^\times}(a_2){\mathbb{I}}_{\varpi^{-1}\o}(a_3)\widehat{\mathbb{I}}_{\o^\times}(a_2a_3\varpi)& = (1-q^{-1}){\mathbb{I}}_{\varpi^{-2}\o^\times}(a_2){\mathbb{I}}_{\varpi\o}(a_3) -q^{-1}{\mathbb{I}}_{\varpi^{-2}\o^\times}(a_2){\mathbb{I}}_{\o^\times}(a_3).
\end{align*}
Direct computation shows that
\begin{align*}
I_1 &= (\beta_2\alpha_3)^{-2}(\alpha_2\beta_3)^2\int_{\F^\times}W_1^\circ({\bf t}(a_3))\mu_3\chi_2(a_3)|a_3|^s{\mathbb{I}}_{\o}(a_3)\,d^\times a_3\\
&=(\beta_2\alpha_3)^{-2}(\alpha_2\beta_3)^2L(s+\tfrac{1}{2},\chi_1\chi_2\mu_3)L(s+\tfrac{1}{2},\mu_1\chi_2\mu_3),\\
I_2 & = q^s(1-q^{-1})(\beta_2\alpha_3)^{-2}(\alpha_2\beta_3)\int_{\F^\times}W_1^\circ({\bf t}(a_3\varpi))\mu_3\chi_2(a_3)|a_3|^s{\mathbb{I}}_{\o}(a_3)\,d^\times a_3-q^{2s-1}(\beta_2\alpha_3)^{-2}\\
&= -(\beta_2\alpha_3)^{-2}\left[(\alpha_1\beta_1)(\alpha_2\beta_3)^2q^{-1}-\alpha_2\beta_3(\alpha_1+\beta_1)q^{s-1/2}+q^{2s-1}\right] L(s+\tfrac{1}{2},\chi_1\chi_2\mu_3)L(s+\tfrac{1}{2},\mu_1\chi_2\mu_3).
\end{align*}
We conclude that
\begin{align*}
(1-q^{-1})\cdot I_1 -(\alpha_1\beta_1)^{-1}\cdot I_2 &= (\beta_2\alpha_3)^{-2}(\alpha_2\beta_3)^{2}(1-\alpha_1^{-1}\alpha_2^{-1}\beta_3^{-1}q^{s-1/2})(1-\beta_1^{-1}\alpha_2^{-1}\beta_3^{-1}q^{s-1/2})\\
&\times L(s+\tfrac{1}{2},\chi_1\chi_2\mu_3)L(s+\tfrac{1}{2},\mu_1\chi_2\mu_3)\\
& = (\beta_2\alpha_3)^{-2}(\alpha_2\beta_3)^{2}\gamma(s+\tfrac{1}{2},\itPi_1\otimes\chi_2\mu_3,\psi)^{-1}.
\end{align*}
This completes the proof of the first assertion.

Now we consider the second assertion.
It is clear that $\itPhi'$ factors through ${\rm Sym}_3(\varpi^{-1}\o) / {\rm Sym}_3(\varpi\o)$.
By Lemma \ref{L:local zeta 2}, we have
\begin{align*}
&Z(W_1^{(1)}\otimes W_2^{(1)}\otimes W_3^{(1)},f_{\itPhi'}^{(s)}) \\
& = q^{-2}(1+q^{-1})^{-2}\gamma(s+\tfrac{1}{2},\itPi_1\otimes\mu_2\mu_3,\psi)^{-1}\\
&\times\int_{\F}dx\int_{(\F^\times)^2}d^\times a_2\,d^\times a_3\,W_1^{(1)}({\bf t}(a_2a_3)J_1{\bf n}(x))
\mu_2\chi_3(a_2)|a_2|^{s}
\mu_3\chi_2(a_3)|a_3|^{s}\widehat{\mathbb{I}}_{\o^\times}(a_2)\widehat{\mathbb{I}}_{\o^\times}(a_3){\mathbb{I}}_{\varpi^{-1}\o}(x).
\end{align*}
Note that $W_1^{(1)}$ is right invariant by ${\bf n}(\varpi^{-1}\o)$. Indeed, for $x \in \varpi^2\o$ and $y \in \varpi^{-1}\o$, we have $1+xy \in \o^\times$ and 
\[
{\bf n}^-(x){\bf n}(y) = {\bf m}(1+xy)^{-1}{\bf n}(y(1+xy)){\bf n}^-(x(1+xy)^{-1}).
\]
Thus 
\[
f_1^{(1)}({\bf n}^-(x){\bf n}(y)) = f_1^{(1)}({\bf n}^-(x(1+xy)^{-1})) = f_1^{(1)}({\bf n}^-(x)).
\]
Therefore, 
\begin{align*}
&\int_{\F}dx\int_{(\F^\times)^2}d^\times a_2\,d^\times a_3\,W_1^{(1)}({\bf t}(a_2a_3)J_1{\bf n}(x))
\mu_2\chi_3(a_2)|a_2|^{s}
\mu_3\chi_2(a_3)|a_3|^{s}\widehat{\mathbb{I}}_{\o^\times}(a_2)\widehat{\mathbb{I}}_{\o^\times}(a_3){\mathbb{I}}_{\varpi^{-1}\o}(x)\\
& = q\int_{(\F^\times)^2}d^\times a_2\,d^\times a_3\,W_1^{(1)}({\bf t}(a_2a_3)J_1)
\mu_2\chi_3(a_2)|a_2|^{s}
\mu_3\chi_2(a_3)|a_3|^{s}\widehat{\mathbb{I}}_{\o^\times}(a_2)\widehat{\mathbb{I}}_{\o^\times}(a_3)\\
& = q^{-1}\int_{(\F^\times)^2}d^\times a_2\,d^\times a_3\,\mu_1\mu_2\chi_3(a_2)|a_2|^{s+1/2}
\mu_1\mu_3\chi_2(a_3)|a_3|^{s+1/2}\widehat{\mathbb{I}}_{\o^\times}(a_2)\widehat{\mathbb{I}}_{\o^\times}(a_3)\mathbb{I}_{\varpi^{-2}\o}(a_1a_2)\\
& = q^{-1}\gamma(s+\tfrac{1}{2},\mu_1\mu_2\chi_3,\psi)^{-1}\gamma(s+\tfrac{1}{2},\mu_1\mu_3\chi_2,\psi)^{-1}.
\end{align*}
This completes the proof.
\end{proof}

\begin{rmk}\label{R:support}
When $2 \nmid q$, it is easy to see that 
\begin{align*}
{\rm supp}(\widehat{\Phi}) \subset \bp 
\o^\times &\varpi\o &  \varpi\o \\
\varpi\o  &\varpi\o & \o^\times \\
\varpi\o  &\o^\times&\varpi\o
\ep,\quad {\rm supp}(\widehat{\Phi}') \subset \bp 
\varpi\o   &\o^\times&  \o^\times \\
\o^\times  &\varpi\o & \o^\times \\
\o^\times  &\o^\times&\varpi\o
\ep.
\end{align*}
In particular, both $\widehat{\Phi}$ and $\widehat{\Phi}'$ are supported in ${\rm Sym}_3(\F)\cap \GL_3(\o)$.
\end{rmk}


\end{document}